\documentclass{article}
\usepackage[utf8]{inputenc}
\usepackage{amsmath,amssymb,amsthm,mathtools}
\usepackage{booktabs,array}
\usepackage{graphicx}
\usepackage[margin=2.6cm]{geometry}
\usepackage{color}
\usepackage[sort]{cite}
\usepackage[colorlinks=true,linkcolor=blue,citecolor=blue,urlcolor=blue]{hyperref}
\usepackage[expansion=false]{microtype}

\newcommand{\R}{\mathbb{R}}
\newcommand{\T}{\mathcal{T}}
\newcommand{\range}{\operatorname{range}}
\newcommand{\supp}{\operatorname{supp}}
\newcommand{\dd}{\,\mathrm{d}}
\newcommand{\adm}{\mathcal{A}}

\newtheorem{theorem}{Theorem}[section]
\newtheorem{proposition}[theorem]{Proposition}
\newtheorem{corollary}[theorem]{Corollary}
\newtheorem{lemma}[theorem]{Lemma}
\theoremstyle{definition}
\newtheorem{definition}[theorem]{Definition}
\theoremstyle{remark}
\newtheorem{remark}[theorem]{Remark}

\title{State-constrained optimal control with the theory of functional connections: squared slacks and junction exponents}

\author{
Emmanuel Tr\'elat\footnote{Sorbonne Universit\'e, Universit\'e Paris Cit\'e, CNRS, Inria, Laboratoire Jacques-Louis Lions, LJLL, F-75005 Paris, France (\texttt{emmanuel.trelat@sorbonne-universite.fr}).}
~and Daniele Mortari\footnote{Emeritus, Aerospace Engineering, Texas A\&M University, College Station TX, 77843-3141, USA (\texttt{mortari@tamu.edu}).}
}
\date{}

\begin{document}
\maketitle

\begin{abstract}
For optimal control problems with an affine state constraint, we combine the theory of functional connections (TFC), which writes every function satisfying prescribed linear conditions as an explicit expression in terms of a free function, with the squared-slack substitution $b - a^\top x(t) = s(t)^2$.
Every trial trajectory then satisfies the endpoint conditions and the path constraint identically.
The accuracy depends on the regularity of the slack function at the junctions with a boundary arc.
Along an optimal trajectory, the constraint residual vanishes like $(t_1 - t)^m$, where the junction exponent $m$ is equal to the order $p$ of the constraint if and only if the $p$-th derivative of the constraint jumps at the junction, and thus a single polynomial slack converges algebraically.
The \emph{multiple-arc TFC method} that we introduce, with unknown junction times and with residuals $(t_1 - t)^m \sigma(t)^2$ on the free arcs, represents the optimal trajectory exactly and converges geometrically when the arcs are analytic.
In a linear-quadratic class, an exact cost identity shows that the cost of every trial trajectory is an upper bound of the optimal cost.
We compare both approaches with direct collocation and indirect shooting on the Bryson-Denham problem, on a rendezvous problem with an approach corridor and on the six-dimensional atmospheric reentry of a space shuttle with three state constraints.
\end{abstract}


\section{Introduction}\label{introduction}

Pure state constraints are common in aerospace problems (altitude and velocity bounds, heating limits during atmospheric flight, approach corridors during rendezvous and landing).
The optimal solutions may contain boundary arcs, along which the constraint is active, joined to free arcs at junction times that are unknown a priori.
At the junctions the costate jumps, and a multiplier measure lives on the boundary arcs \cite{BonnansHermant2009AIHP,Bonnard2003Reentry,BrysonDenhamDreyfus1963,HartlSethiVickson1995,JacobsonLeleSpeyer1971,Maurer1977}.

Direct methods discretize the trajectory and impose the constraint at the nodes.
They usually converge from crude initial guesses, but their accuracy is limited by the order of the discretization, the constraint is not controlled between the nodes and the discrete optimal value is not a bound of the true one \cite{Betts1998,Betts2010,ByczkowskiRao2024,HargravesParis1987,Kelly2017,Trelat2012JOTA}.
Indirect methods solve the boundary value problem of the Pontryagin maximum principle, junction conditions included.
They reach the precision of the integrator, but they require the arc structure and a good initial costate, which homotopy methods and initializations from a direct solution help to provide \cite{BonnansHermant2007SICON,BonnardFaubourgTrelat2005,BulirschMontronePesch1991,CaillauCotsGergaud2012,CaillauCotsMartinon2022ct,HaberkornMartinonGergaud2004,MaurerGillessen1975,Pesch1994,Trelat2012JOTA}.

\medskip
The theory of functional connections (TFC) provides an explicit expression of all functions satisfying a family of linear constraints, in the form of a free function corrected by switching functions \cite{LeakeJohnstonMortari2022,Mortari2017Connections,MortariLeake2019}.
The free function is expanded on a polynomial basis and determined by least squares, so that the constraints are satisfied identically and only the remaining equations are solved approximately \cite{Mortari2017LS,MortariJohnstonSmith2019,Schiassi2021XTFC}.
This construction has been applied to the boundary value problems of the maximum principle arising in orbit transfers, relative motion, landing and reentry \cite{DAmbrosio2022Landing,DrozdFurfaro2025,DrozdFurfaroDAmbrosio2024,Drozd2021,FurfaroMortari2020,Johnston2020Descent}.
Inequality constraints have been treated in TFC by sigmoid formulas with a quantified overrun \cite{JohnstonLeakeMortari2020Ineq}, or by Heaviside switches and clipping \cite{LeakeThesis2021,Mortari2025MinParam}.

The tools that we use are classical.
The squared-slack substitution \cite{JacobsonLele1969,Valentine1937} is standard in nonlinear programming \cite{Bertsekas1999,DingWright2025,Taylor1973}, and it has been used with the adjoint method in optimal control \cite{EichmeirNachbagauerSteiner2025}.
We also use the piecewise constrained expressions of \cite{deAlmeida2026Segmentation,JohnstonMortari2021Hybrid}, and the restriction of the cost to a family of admissible trajectories, which is the Ritz method of the calculus of variations.

Let us mention other constructions dealing with the feasibility of the trajectory on the whole horizon.
Saturation functions parametrize only the interior of the constraint set \cite{GraichenPetit2009}, and B-spline parametrizations guarantee the constraint through sufficient conditions on the coefficients, with some conservatism \cite{VanLoock2015}.
Flatness-based trajectory generation usually imposes the inequalities at collocation points \cite{Fliess1995,MilamMushambiMurray2000}.
Interior-point methods in function space approach the constraint from inside, and convergence results of the primal and dual variables have been obtained for them under state constraints \cite{Malisani2023SICON,Malisani2024COCV}.

\medskip
In this paper, we study the combination of the squared slack with the constrained expressions of TFC, for a constraint that is affine in the state, and the consequences of the regularity of the slack at the junctions.
We consider a linear inequality constraint $a^\top x(t) \leq b$, because the substitution $b - a^\top x(t) = s(t)^2$ then turns it into a linear equality, to which TFC applies in coordinates adapted to the constraint.
In these coordinates, the endpoint conditions on the state become endpoint conditions on the slack, and every trial trajectory satisfies the endpoint conditions and the path constraint identically (an active constraint corresponds to a zero of the slack).
The admissibility of all trial trajectories and the cost bounds are guaranteed in the linear-quadratic class, where the dynamics are eliminated, whereas the reentry example tests the representation of the constraints and the parametrization by unknown junction times without these guarantees.

In nonlinear programming, the slack is a finite-dimensional variable and only its value matters.
In optimal control, it becomes a function of time, whose derivatives enter the control and the cost, and its regularity determines what can be represented by a polynomial expansion.

On a free arc adjacent to a boundary arc entered at $t_1$, the residual $y^\star(t) = b - a^\top x^\star(t)$ of the optimal trajectory $x^\star(\cdot)$ vanishes like $(t_1 - t)^m$.
We show that the junction exponent $m$ is at least the order $p$ of the constraint along the control system, with equality if and only if the $p$-th derivative of the constraint along the trajectory jumps at the junction.
Hence the slack $\sqrt{y^\star}$ behaves like $(t_1 - t)^{m/2}$, and it is H\"older of exponent $m/2$ and no better.
Now, a slack expanded on a finite family of analytic functions on the whole horizon, as in the usual TFC computations, cannot vanish identically on the boundary arc, and the best uniform approximation of the optimal slack by polynomials of degree $N$ has an error of the order of $N^{-m/2}$.
We call this approach the \emph{single-arc TFC method}, by analogy with single shooting, in which one integration covers the whole horizon.
Its obstructions (in representation, accuracy and conditioning) are gathered in Remark~\ref{single-obstructions}.

As a remedy, we introduce a \emph{multiple-arc TFC method}, named by analogy with multiple shooting (which replaces the single integration by several integrations glued at nodes).
Here, the horizon is cut at the junction times, which become unknowns, the boundary arcs are exact, and the residual on each free arc is $(t_1 - t)^m \sigma(t)^2$, with a factor $\sigma(\cdot)$ that is as regular as the free arc itself.
The family of trajectories generated by this method contains the optimal trajectory, and the method converges geometrically when the arcs are analytic.
For fully actuated linear systems with quadratic cost, the dynamics are eliminated and every trial trajectory is admissible, and thus its cost is an upper bound of the optimal cost.
Moreover, an exact identity involving the multiplier measure of the constraint expresses the cost error, and explains why the residual must be kept nonnegative on the free arcs (by the square), and not only tangent at the junctions.

\medskip
As for multiple shooting, the multiple-arc TFC method requires the number of boundary arcs, their order and the junction exponents to be known before the junction times can be taken as unknowns.
Since this structure is unknown a priori, we proceed in two steps.
In a first step, we run the single-arc TFC method, which needs no structure.
A threshold on its slack function locates the boundary arcs roughly, and gives their number and approximate junction times.
In a second step, we initialize the multiple-arc TFC method with these values, expecting that this rough approximation lies in its basin of convergence (this was the case in all linear-quadratic examples of this paper).
We also describe how the exponents are chosen among two candidates, and why additional cuts inside an arc are harmless.
When the dynamics are not eliminated, the first step can still be performed in a residual form, and the structure and the initialization can also be obtained from a direct method or from a geometric analysis of the extremals.

We study three examples in detail.
The Bryson-Denham problem \cite{BrysonDenhamDreyfus1963,BrysonHo1975} has a closed-form solution with three regimes.
We also consider a forced variant, which keeps the closed form but removes the accidental exactness of piecewise-cubic schemes.
The rendezvous problem in the Clohessy-Wiltshire frame with an affine approach corridor \cite{ClohessyWiltshire1960,Drozd2021} has a boundary arc along which the optimal motion is nontrivial.
The atmospheric reentry of a space shuttle \cite{Bonnard2003Reentry,BonnardFaubourgTrelat2005,BonnardTrelat2002COCV,Trelat2003DCDS} has six states, a bank-angle control, a free final time and three nonlinear state constraints, on the thermal flux, the normal acceleration and the dynamic pressure.
This example shows what remains of the construction outside the linear-quadratic class: the exact representation of each constraint on the arcs adjacent to its boundary arcs (by a nonlinear pullback), the multiple-arc TFC method with unknown junction times, and the single-arc TFC method in a residual form, used to find the structure.
But the trial trajectories are not admissible anymore, and the bounds are lost.
We compare the TFC methods with direct collocation, solved by \texttt{IpOpt} \cite{WachterBiegler2006} through \texttt{CasADi} \cite{Andersson2019CasADi}, and with indirect shooting on the junction conditions.

\medskip
The paper is organized as follows.
Section~\ref{expressions} is devoted to constructing the two methods and to proving the result on the junction exponent.
In Section~\ref{analysis}, we establish the cost identity and derive its consequences.
The numerical examples are presented in Section~\ref{examples}, and Section~\ref{conclusion} concludes the paper.

\section{Squared slacks and the multiple-arc TFC method}\label{expressions}

\subsection{Linear constraints and squared slacks}\label{finite}

Let $X$ and $Y$ be vector spaces, and let $L : X \to Y$ be a linear map with a linear right inverse $R : Y \to X$, i.e., $LR = I_Y$.
Given $c \in Y$, the constrained expression of TFC is defined by
\begin{equation}\label{abstract-tfc}
\T_c[g] = g + R(c - Lg) \qquad \forall g \in X .
\end{equation}
Since $L\T_c[g] = c$ for every $g$ and $\T_c[u] = u$ whenever $Lu = c$, the range of $\T_c$ is the affine set $\{u \in X \mid Lu = c\}$.
When $Y = \R^q$ and $L = (L_1, \dots, L_q)$, a right inverse is built from support functions $\phi_1, \dots, \phi_q \in X$ such that the matrix $M_{ij} = L_i \phi_j$ is invertible, and \eqref{abstract-tfc} takes the switching form $\T_c[g] = g + \sum_i \psi_i (c_i - L_i g)$, with $L_j \psi_i = \delta_{ij}$.
Note that only the $q \times q$ matrix $M$ is inverted and that no basis of $\ker L$ is needed: this is the practical advantage of the construction in function spaces \cite{MaiMortari2022,MortariLeake2019}.

Let us now consider an inequality $Ax \leq b$ in $\R^n$, with $A \in \R^{m \times n}$ of full row rank, and let $A^\dagger = A^\top (AA^\top)^{-1}$.
The squared-slack substitution $Ax + s^2 = b$, with $s \in \R^m$ and $s^2$ taken componentwise, turns the inequality into a linear equality with a parameter-dependent target, and by \eqref{abstract-tfc} we get
$$
\T_b[g, s] = A^\dagger (b - s^2) + (I_n - A^\dagger A) g, \qquad \range \T_b = \{x \in \R^n \mid Ax \leq b\} :
$$
every pair $(g, s)$ provides a feasible point, and every feasible point is reached, with $s_i = \pm \sqrt{b_i - (Ax)_i}$.
An active constraint is reached at a finite parameter value ($s_i = 0$), unlike with exponential slacks, which parametrize only the open set, and without the nonsmoothness of Heaviside or clipping constructions.

But this parametrization is singular on the boundary.
Indeed, in the reduced coordinates $x(s, z) = A^\dagger (b - s^2) + Qz$, with $Q$ a basis of $\ker A$, the $i$-th column of $D_s x = -2A^\dagger \operatorname{diag}(s)$ vanishes when $s_i = 0$.
Actually, every smooth parametrization of a closed inequality set, defined on an open parameter domain, has this property, since $h \circ \Psi$ has a maximum at any interior parameter where $h(\Psi(\bar z)) = 0$.
The active normal direction is lost at first order, and the sign of the multiplier is recovered at second order.

\begin{lemma}[second-order stationarity and KKT]\label{kkt-transfer}
Let $f \in C^2(\R^n)$ and $F(s, z) = f(x(s, z))$.
If $(s^\star, z^\star)$ is a stationary point of $F$ at which $\nabla^2 F$ is positive semidefinite, then $x^\star = x(s^\star, z^\star)$ satisfies the KKT conditions of $\min \{f(x) \mid Ax \leq b\}$, with the multiplier $\lambda^\star$ defined by $\nabla f(x^\star) + A^\top \lambda^\star = 0$.
\end{lemma}

\begin{proof}
Stationarity in $z$ gives $Q^\top \nabla f(x^\star) = 0$, hence a unique $\lambda^\star$.
Stationarity in $s_i$ gives $2 s_i^\star \lambda_i^\star = 0$, which is complementary slackness since $(Ax^\star - b)_i = -(s_i^\star)^2$.
On an active index, $\partial^2_{s_i s_i} F = -2 (A^\dagger e_i)^\top \nabla f(x^\star) = 2\lambda_i^\star \geq 0$.
\end{proof}

Under strict complementarity, the curvature $2\lambda_i^\star$ is positive in every active slack direction, and it vanishes with the multiplier.
These properties of squared-variable formulations are classical \cite{Bertsekas1999,DingWright2025,FukudaFukushima2017,Taylor1973}.
Lemma~\ref{kkt-transfer} answers, at second order and in finite dimension, the observation attributed to Polak in \cite{JacobsonLele1969}, that a stationary point of the transformed problem need not be stationary for the original one.

\subsection{Adapted coordinates and the single-arc TFC method}\label{adapted}

Let $I = [0, T]$ and let $x(\cdot) : I \to \R^n$ be a trajectory of class $C^k$.
Two kinds of constraints act on $x(\cdot)$.
The first ones are linear endpoint conditions $Lx = c$, which prescribe values and consecutive derivatives of $x(\cdot)$ at $t = 0$ and $t = T$.
For them, the support functions can be taken to be the monomials $1, t, \dots, t^{q-1}$, and a right inverse exists, because the matrix $M_{ij} = L_i \phi_j$ of Section~\ref{finite} is then the matrix of Hermite interpolation at $0$ and $T$ (a confluent Vandermonde matrix), which is invertible since $0 \neq T$.
The second one is an affine path constraint
$$
a^\top x(t) \leq b \qquad \forall t \in I, \qquad a \in \R^n \setminus \{0\}, \ b \in \R .
$$

Our objective is to vary a set of coefficients freely, without having to check the endpoint conditions and the pointwise inequality at each step.
For an affine constraint, the inversion is global and explicit.
Let $Q \in \R^{n \times (n-1)}$ have orthonormal columns spanning $a^\perp$.
The affine change of coordinates
\begin{equation}\label{adapted-coordinates}
y(t) = b - a^\top x(t) \in \R, \qquad z(t) = Q^\top x(t) \in \R^{n-1}, \qquad x(t) = \frac{a}{\vert a \vert^2} (b - y(t)) + Qz(t) ,
\end{equation}
is a bijection of $C^k(I; \R^n)$ onto $C^k(I; \R) \times C^k(I; \R^{n-1})$, with $y^{(j)} = -a^\top x^{(j)}$ and $z^{(j)} = Q^\top x^{(j)}$ for $j \geq 1$, and the path constraint becomes $y(t) \geq 0$.
We assume that the endpoint conditions separate in these coordinates, i.e., $L = (L_y, L_z)$ and $c = (c_y, c_z)$, which is the case when full state vectors (or derivatives thereof) are prescribed.

The substitution $y(t) = s(t)^2$ imposes the sign of the residual, and a constrained expression for $s(\cdot)$ imposes the endpoint data.
If only the values $s(0) = s_0$ and $s(T) = s_T$ are prescribed, the constrained expression is the linear correction
$$
s(t) = g(t) + \Big(1 - \frac{t}{T}\Big) \big(s_0 - g(0)\big) + \frac{t}{T} \big(s_T - g(T)\big) ,
$$
valid for every free function $g(\cdot)$, and we have $y(t) = s(t)^2 \geq 0$ on $I$.
When derivatives are prescribed as well, the correction becomes the Hermite interpolant of the prescribed values and derivatives.

Let us now determine which data on $s(\cdot)$ encode the data $L_y y = c_y$ on $y(\cdot)$.
At an endpoint $t_i$ where $y(t_i) = y_i > 0$ is prescribed, the slack satisfies $s(t_i) = \epsilon_i \sqrt{y_i}$, with a sign $\epsilon_i \in \{-1, 1\}$.
Once the sign is fixed, the higher-order conditions transfer uniquely, since $s \mapsto s^2$ is a local diffeomorphism of jets where $s \neq 0$ (for instance, $\dot s(t_i) = \dot y(t_i) / (2 \epsilon_i \sqrt{y_i})$).

If the constraint is active at an endpoint, we only consider the case of a second-order constraint with the data $y(t_i) = 0$ and $\dot y(t_i) = 0$ (zero normal velocity on the constraint boundary), as for the target of the rendezvous problem of Section~\ref{rendezvous}.
Since a $C^1$ slack with $s(t_i) = 0$ satisfies both conditions, the transferred condition is only $s(t_i) = 0$, with no sign attached.
Data with $y(t_i) = 0$ and a nonzero derivative are not compatible with a $C^1$ slack, and we do not consider them.
We call \emph{transferred data} the resulting conditions $\tilde L_y s = \tilde c_y$.
They depend on the endpoint signs $\epsilon = (\epsilon_0, \epsilon_T)$, at the endpoints where $y$ is prescribed positive.

\begin{definition}[single-arc TFC method]\label{single}
Let $\T_{\tilde c_y}$ and $\T_{c_z}$ be constrained expressions on $I$ for the transferred data on $s(\cdot)$ and for the data on $z(\cdot)$.
The \emph{single-arc TFC method} uses one single constrained expression on the whole horizon of time,
$$
\mathcal{X}[g_s, g_z] = \frac{a}{\vert a \vert^2} \big(b - \T_{\tilde c_y}[g_s]^2\big) + Q \T_{c_z}[g_z], \qquad (g_s, g_z) \in C^k(I) \times C^k(I; \R^{n-1}) .
$$
\end{definition}

The name is chosen by analogy with single shooting, in which one integration covers the whole horizon and the unknowns are the initial data.
Here, one constrained expression covers the whole horizon, and the unknowns are the coefficients of its free functions.
The multiple-arc TFC method of Section~\ref{multiple} cuts the horizon into arcs glued by junction values, as multiple shooting glues its integrations by matching conditions at the nodes.

Let $\adm$ be the set of trajectories $x(\cdot) \in C^k(I; \R^n)$ such that $Lx = c$ and $a^\top x(t) \leq b$ on $I$.
The following result characterizes the elements of $\adm$ that are reached by the single-arc TFC method, in terms of the regularity of the slack function and of the endpoint signs.

\begin{proposition}[exact range]\label{exact-range}
\begin{enumerate}
\item[(i)] One has
$$
\range \mathcal{X} = \{x(\cdot) \in \adm \ \mid\ b - a^\top x(t) = s(t)^2 \text{ on } I \text{ for some } s(\cdot) \in C^k(I) \text{ with } \tilde L_y s = \tilde c_y\} .
$$
In particular, every trial trajectory satisfies the endpoint conditions and the path constraint identically, on the whole interval.
\item[(ii)] Let $x(\cdot) \in \adm$ with $b - a^\top x(t) > 0$ on $I$.
Then $x(\cdot) \in \range \mathcal{X}$ if and only if $\epsilon_0 \epsilon_T = 1$, and then $s = \epsilon_0 \sqrt{b - a^\top x}$.
\item[(iii)] Let $x(\cdot) \in \adm$ with $y = b - a^\top x$ positive at both endpoints and with finitely many zeros $\tau_j$ in $(0, T)$, all of even order $2r_j$ with $r_j \leq k$ and $y / (t - \tau_j)^{2r_j}$ of class $C^k$ and positive near $\tau_j$.
Then $x(\cdot) \in \range \mathcal{X}$ if and only if $\epsilon_0 \epsilon_T = (-1)^\nu$, where $\nu$ is the number of zeros $\tau_j$ with $r_j$ odd.
\end{enumerate}
\end{proposition}

\begin{proof}
(i) If $x = \mathcal{X}[g_s, g_z]$, then $L_z z = c_z$ by the retraction property of $\T_{c_z}$, $y = s^2 \geq 0$, and $L_y y = c_y$ because the transferred data were built so that $L_y(s^2) = c_y$ whenever $\tilde L_y s = \tilde c_y$.
Conversely, if $x(\cdot) \in \adm$ and $y = s^2$ with $\tilde L_y s = \tilde c_y$, then $\T_{\tilde c_y}[s] = s$ and $\T_{c_z}[z] = z$, hence $\mathcal{X}[s, z] = x$.\\
(ii) If $y > 0$ then a $C^k$ function $s(\cdot)$ with $s^2 = y$ is $\pm \sqrt{y}$ with a constant sign; the transferred data impose $s(0) = \epsilon_0 \sqrt{y(0)}$ and $s(T) = \epsilon_T \sqrt{y(T)}$, compatible with a constant sign if and only if $\epsilon_0 \epsilon_T = 1$, and $s = \epsilon_0 \sqrt{y}$ satisfies the transferred data because the jets of $s$ at the endpoints are determined by those of $y$ and by the sign.\\
(iii) Near $\tau_j$ the $C^k$ square roots of $y$ are $\pm (t - \tau_j)^{r_j} \sqrt{y / (t - \tau_j)^{2r_j}}$, and the sign of $s$ changes across $\tau_j$ if and only if $r_j$ is odd; a $C^k$ function $s(\cdot)$ with $s^2 = y$ and the prescribed endpoint signs exists if and only if the number of sign changes has the parity of $\epsilon_0 \epsilon_T$.
\end{proof}

Proposition~\ref{exact-range} shows that a single $C^k$ slack function $s(\cdot)$ represents the trajectories lying strictly inside the constraint, and the trajectories touching the boundary at isolated points of even order, provided that the endpoint signs have the right parity.
If an endpoint is active, then no sign is prescribed there and the parity conditions of (ii) and (iii) disappear, but the regularity requirement on the root remains.

\begin{remark}[obstructions to the single-arc TFC method]\label{single-obstructions}
The single-arc TFC method meets its limits on a boundary arc, i.e., on a subinterval on which $y \equiv 0$, and we gather these limits below.
The multiple-arc TFC method of Section~\ref{multiple} removes the first four of them, but it requires the structure of the solution.
\begin{enumerate}
\item[(i)] Regularity.
Writing $y = s^2$ with $s(\cdot) \in C^k$ requires $\sqrt{y}$ to be of class $C^k$ across the junctions.
But, at a junction of an optimal trajectory, $\sqrt{y^\star}$ is $C^{q, \alpha}$ and no better, with $q = \lceil m/2 \rceil - 1$ (see Proposition~\ref{junction-exponent} below).
Hence the optimal trajectory belongs to the range only for $k \leq q$, i.e., for $k \leq 1$ for the Bryson-Denham problem of Section~\ref{bryson}.
\item[(ii)] Representation.
Smooth slack functions vanishing on an interval do exist, but all their derivatives vanish at the junctions, whereas the slack of an optimal trajectory has a finite-order contact.
Moreover, a finite sum of analytic functions, which is what TFC computations use (polynomials and Fourier modes in \cite{LeakeJohnstonMortari2022,MortariLeake2019}, networks with analytic activations in \cite{LeakeMortari2020Deep,Schiassi2021XTFC}), cannot vanish on an interval without vanishing everywhere.
With such a free function, the optimal trajectory lies outside the range whatever the degree may be.
The slack can only be small on the boundary arc, and its best polynomial approximation converges at the algebraic rate $N^{-m/2}$ (Corollary~\ref{jackson}).
\item[(iii)] Cost.
In the linear-quadratic class, the error identity of Proposition~\ref{ritz-theorem} contains the term $\int (b - a^\top \xi) \dd\mu$, which does not vanish when the slack is small but nonzero where the multiplier measure lives.
The cost error is then algebraic (we observe a decay like $N^{-3}$ when $m = 3$), and the control error, which involves $\ddot s$, decays only like $N^{-1}$ (Section~\ref{bryson}).
\item[(iv)] Conditioning.
The derivative of the control with respect to the slack vanishes with $s$ and its derivatives, and thus the Gauss-Newton iteration degenerates once the boundary arc is resolved (Section~\ref{examples}).
This is the functional counterpart of the singularity of Section~\ref{finite}.
For the same reason, the junction times are located only to a few hundredths of the horizon, where the slack becomes small.
\item[(v)] Parity.
The endpoint signs are a discrete choice.
A wrong choice does not produce infeasible trajectories, but it restricts the ansatz to a family that cannot contain the optimum, and Remark~\ref{parity-gap} quantifies the resulting gap in cost.
On a boundary arc, the slack vanishes identically, so that its signs on the two adjacent free arcs are independent and both parities are admissible.
\end{enumerate}
A free function that is analytic by pieces, with its breakpoints at the junctions and the right vanishing order there, has none of the defects (i) to (iv), which motivates the multiple-arc construction of Section~\ref{multiple}.
\end{remark}

\subsection{The junction exponent}\label{junction}

Let us consider a control system $\dot x(t) = f(x(t), u(t))$, with $f$ of class $C^p$, a state constraint $c(x) \leq 0$, with $c$ of class $C^p$, and an admissible pair $(x^\star(\cdot), u^\star(\cdot))$, with $u^\star(\cdot)$ piecewise continuous.
Recall that the constraint is said to be of order $p$ along the control system if the iterated derivatives $c^{(j)}(x) = \nabla c^{(j-1)}(x)^\top f(x, u)$ do not depend on $u$ for $j < p$ and $c^{(p)}(x, u)$ does \cite{Bonnard2003Reentry,BrysonDenhamDreyfus1963,HartlSethiVickson1995,Maurer1977}.
A boundary arc is an interval $[t_1, t_2]$ on which $c(x^\star(t)) \equiv 0$, so that $c^{(p)}(x^\star(t), u^\star(t)) \equiv 0$ there.

We set $y^\star(t) = -c(x^\star(t)) \geq 0$ and say that the entry junction $t_1$ has a \emph{finite-order contact} if, for some $\delta > 0$,
\begin{equation}\label{finite-contact}
y^\star(t) = (t_1 - t)^m \varphi(t) \quad \text{on } [t_1 - \delta, t_1], \qquad m \in \mathbb{N}^\ast, \quad \varphi \text{ continuous}, \quad \varphi(t_1) > 0 .
\end{equation}
The integer $m$, which is then unique, is called the \emph{junction exponent}.
Note that the contact of an analytic free arc is of finite order, but that $C^\infty$ regularity does not imply it: for instance, $y^\star(t) = \exp(-1/(t_1 - t)^2)$ on the free arc, extended by zero, is $C^\infty$ and admits no such $m$.

\begin{proposition}[junction exponent]\label{junction-exponent}
Let $[t_1, t_2]$ be a boundary arc entered at $t_1$ with a finite-order contact \eqref{finite-contact}.
Then:
\begin{enumerate}
\item[(i)] $m \geq p$, and $m = p$ if and only if $c^{(p)}\big(x^\star(t_1), u^\star(t_1^-)\big) \neq 0$, that is, if and only if the $p$-th derivative of the constraint along the trajectory jumps at the junction; in particular $m \geq p + 1$ whenever $u^\star(\cdot)$ is continuous at $t_1$;
\item[(ii)] let $q = \lceil m/2 \rceil - 1$ and $\alpha = m/2 - q$; if $\varphi \in C^{q+1}$ near $t_1$, the slack $s^\star = \sqrt{y^\star}$ on the free arc, extended by $0$ on the boundary arc, belongs to $C^{q, \alpha}$ at $t_1$, and to no better H\"older class: for $m$ odd it is not $C^{q, \alpha'}$ for $\alpha' > \alpha = 1/2$, and for $m$ even it is $C^{q, 1}$ but not $C^{q+1}$;
\item[(iii)] if $\varphi \in C^k$ near $t_1$, then $\sigma^\star = \sqrt{\varphi}$ is of class $C^k$ and positive there.
The same holds at the exit junction.
\end{enumerate}
\end{proposition}

\begin{proof}
For $j < p$ the function $t \mapsto c^{(j)}(x^\star(t))$ depends on the state only, is continuously differentiable, and vanishes on the boundary arc; on the free arc it equals $-(y^\star)^{(j)}(t)$.
Hence $y^\star$ is of class $C^p$ on the closed free arc, with $(y^\star)^{(j)}(t_1^-) = 0$ for $j < p$, and $(y^\star)^{(p)}(t_1^-) = -c^{(p)}(x^\star(t_1), u^\star(t_1^-))$ since $u^\star$ has a left limit.
By Taylor's formula, we have $y^\star(t) = (y^\star)^{(p)}(t_1^-) (t - t_1)^p / p! + \mathrm{o}((t_1 - t)^p)$, and comparing with \eqref{finite-contact}, we get $m \geq p$, with $m = p$ exactly when $(y^\star)^{(p)}(t_1^-) \neq 0$.
On the arc side, $c^{(p)}(x^\star(t_1), u^\star(t_1^+)) = 0$, and thus the continuity of $u^\star$ forces $m \geq p + 1$, which proves (i).
For (ii), we have $s^\star = (t_1 - t)_+^{m/2} \sqrt{\varphi}$, with $\sqrt{\varphi} \in C^{q+1}$ positive.
The function $(t_1 - t)_+^{m/2}$ is of class $C^{q, \alpha}$ at $t_1$, its derivative of order $q$ behaving like $(t_1 - t)_+^\alpha$, and we get the H\"older class and its sharpness.
Statement (iii) follows from the positivity of $\varphi$.
\end{proof}

The converse of the last assertion of (i) is false: a jump of the control in a direction that $c^{(p)}$ does not see leaves $m > p$.
For the Bryson-Denham and rendezvous problems of Section~\ref{examples}, we have $p = 2$, the control is continuous at the junction and the derivative of its normal component is not, and thus $m = 3$ and the slack is $C^{1, 1/2}$.
For a double integrator with a velocity bound, we have $p = 1$.
In minimum time, the control jumps from a bang value to the boundary control $0$, hence $m = 1$ and the slack is H\"older of exponent $1/2$.
In minimum energy, the Hamiltonian forces $u^\star(t_1^-) = 0$, so that $m = 2$ and the slack is Lipschitz with a kink.
For the reentry problem of Section~\ref{reentry}, the bank angle jumps at every junction in a direction that is seen by the second derivative of the constraint, and thus $m = p = 2$.

A touch point of even order $2r$ is different: there, $y^\star$ vanishes on both sides, and a sign-changing slack $(t - \tau)^r \sqrt{y^\star / (t - \tau)^{2r}}$ can be chosen.
For instance, at the touch point of the Bryson-Denham problem with $1/6 < \ell < 1/4$, setting $h = t - 1/2$, we have $y^\star = (12\ell - 2) h^2 + (4 - 16\ell) \vert h \vert^3$.
The quotient $y^\star / h^2$ is Lipschitz and not $C^1$, because the derivative of the control jumps, its nonnegative root has a cusp, and the signed root $h \sqrt{y^\star / h^2}$ is $C^{1,1}$.
At $\ell = 1/6$, the quadratic term vanishes, the contact is cubic on both sides and the signed root is only $C^{1, 1/2}$.

The rate of polynomial approximation of such a slack is given by the Jackson and Bernstein theorems \cite[Chapter~7]{DeVoreLorentz1993}, \cite[Chapter~7]{Trefethen2019ATAP}.

\begin{corollary}[polynomial approximation of the slack]\label{jackson}
In the setting of Proposition~\ref{junction-exponent}(ii), let $K$ be a compact interval containing $t_1$ in its interior on which $s^\star$ has no other singularity than the junction.
The error of best uniform approximation of $s^\star$ on $K$ by polynomials of degree $N$ satisfies $E_N(s^\star) \leq C N^{-m/2}$, and this rate is sharp for the model singularity $(t_1 - t)_+^{m/2}$.
\end{corollary}

Note that this rate is a property of polynomial approximation on a fixed interval, and not of the singularity alone.
No finite family of analytic functions can represent the boundary arc, and a piecewise-polynomial basis can do it only if a knot sits at the junction, with the right vanishing order.
The multiple-arc TFC method of the next section places such knots, the junction times being unknowns.

\begin{remark}[the transformation technique of \cite{JacobsonLele1969}]
The transformation technique of \cite{JacobsonLele1969} writes the constraint as $c(x) + \frac{1}{2} \alpha^2 = 0$, differentiates it $p$ times, and takes $\alpha^{(p)}$ as the new control and $\alpha, \dots, \alpha^{(p-1)}$ as new states.
Every control of the transformed problem then yields a feasible trajectory, and the boundary arcs become singular arcs.
Proposition~\ref{junction-exponent} also describes the slack of this technique.
Indeed, near a junction, $\alpha = \sqrt{2 y^\star}$ behaves like $(t_1 - t)^{m/2}$, and thus the new control behaves like $(t_1 - t)^{m/2 - p}$.
For the Bryson-Denham problem ($m = 3$, $p = 2$), the new control has an inverse square-root singularity at the junction, which is integrable, but unbounded and not square integrable.
Therefore a transformed problem posed with a bounded new control excludes this junction, and can only approach it.
Note that the examples of \cite{JacobsonLele1969} approach the constraint boundary without reaching it, so that the singularity does not appear there.
Our remark concerns the function space of the new control, and not the technique itself.
\end{remark}

\subsection{The multiple-arc TFC method}\label{multiple}

Let us assume that the optimal trajectory contains a boundary arc $[t_1, t_2]$, with $0 < t_1 < t_2 < T$ or $t_2 = T$, and with junction exponents $m_1$ at the entry and $m_2$ at the exit, i.e., on the free arcs, $y^\star$ vanishes at the junctions like $(t_1 - t)^{m_1}$ and $(t - t_2)^{m_2}$.
We split the horizon at the junction times into arcs, each of them carrying its own constrained expressions and free functions, and we glue the arcs by junction values, which are unknowns of the reduced problem.
The junction times themselves are unknowns, and the residual on each free arc carries the vanishing order of its junction.
By analogy with multiple shooting, we call the result the \emph{multiple-arc TFC method}.
The junction values play the role of the matching conditions, with the difference that the nodes sit at the junctions with the boundary arcs and are unknowns.

Gluing constrained expressions along a partition of the horizon is not new.
Piecewise constrained expressions were introduced in \cite{JohnstonMortari2021Hybrid}, with the junctions at the known switching times of a hybrid system, and, under the name of segmentation, in \cite{deAlmeida2026Segmentation}, where a uniform subdivision of the horizon with position and velocity matching is used to improve accuracy, efficiency and robustness.
The novelty here is the use of the partition to represent the boundary arcs exactly, and it relies on three features that do not appear in the above-mentioned works.
First, the cuts sit at the junctions with the boundary arcs, and not on a prescribed grid.
Second, the junction times are unknowns of the reduced problem, determined together with the coefficients, and not given event times.
Third, the constrained expressions on the free arcs carry the factor $(t_1 - t)^{m_1}$ of the junction exponent, so that every trial trajectory is feasible on the whole horizon and the optimal trajectory belongs to the family (Proposition~\ref{multiple-range}), which no piecewise expression on a prescribed grid achieves.

\begin{definition}[multiple-arc TFC method]\label{multiple-def}
Let $\theta = (t_1, t_2, \zeta)$, where $\zeta$ collects the junction values $z^{(j)}(t_1)$, $z^{(j)}(t_2)$ required by the continuity conditions.
On the free arc $[0, t_1]$ set
\begin{equation}\label{F2-free-arc}
y(t) = (t_1 - t)^{m_1} \sigma_1(t)^2, \qquad \sigma_1 = \T[g_1], \qquad z = \T[g_{z,1}] ,
\end{equation}
where the constrained expression of $\sigma_1(\cdot)$ carries the endpoint data at $t = 0$ transferred through the factor $(t_1 - t)^{m_1}$ and that of $z(\cdot)$ carries the endpoint data at $0$ and the junction values at $t_1$; on the boundary arc set $y \equiv 0$ and $z = \T[g_b]$ with the junction values at $t_1$ and $t_2$; on $[t_2, T]$, when $t_2 < T$, proceed symmetrically with $y(t) = (t - t_2)^{m_2} \sigma_2(t)^2$.
The trajectory $x = \mathcal{X}^{(2)}_\theta[g]$ is recovered from \eqref{adapted-coordinates} on each arc.
\end{definition}

Note that the factor $(t_1 - t)^{m_1}$ contains the tangency conditions of the junction ($y^{(j)}(t_1^-) = 0$ for $j < m_1$), whatever $\sigma_1(\cdot)$ may be, and that the boundary arc is exact.
When $n \geq 2$, the constraint manifold has internal degrees of freedom, and the free function $g_b$ of $z(\cdot)$ on the boundary arc is an unknown.
In other words, the motion along the constraint is optimized like on any other arc, as in the rendezvous problem.
When $n = 1$, the boundary arc is entirely determined.

\begin{proposition}[range of the multiple-arc TFC method]\label{multiple-range}
For fixed $\theta$, the range of $\mathcal{X}^{(2)}_\theta$ is the set of trajectories, of class $C^k$ on each arc of the partition $\{0, t_1, t_2, T\}$, satisfying the endpoint conditions, the continuity conditions encoded in $\zeta$, $y \equiv 0$ on $[t_1, t_2]$, and such that $y / (t_1 - t)^{m_1}$ and $y / (t - t_2)^{m_2}$ are squares of $C^k$ functions on $[0, t_1]$ and $[t_2, T]$ satisfying the transferred data.
Every element of the range satisfies the path constraint on $I$, and its normal coordinate $y$ is of class $C^{\min(k, m_i - 1)}$ across the junction $t_i$.
If $\zeta$ encodes the continuity of $z$ and of its derivatives up to order $r - 1$, if $k \geq r$ and $m_1, m_2 \geq r$, the range is contained in $C^{r-1}(I; \R^n) \cap H^r(I; \R^n)$.
If $x^\star(\cdot) \in \adm$ has a boundary arc $[t_1^\star, t_2^\star]$ with junction exponents $m_1, m_2$ and if $y^\star / (t_1^\star - t)^{m_1}$ and $y^\star / (t - t_2^\star)^{m_2}$ are positive and of class $C^k$ up to the junctions, then $x^\star(\cdot)$ belongs to the range for $\theta^\star = (t_1^\star, t_2^\star, \zeta^\star)$.
\end{proposition}

The proof is similar to that of Proposition~\ref{exact-range}, arc by arc.
The regularity of $y$ across $t_i$ follows from $y^{(j)}(t_i^\mp) = 0$ for $j < m_i$ and from $y \equiv 0$ on the boundary arc, and the tangential coordinate is $C^{r-1}$ across the junctions and $C^k$ on each arc, hence in $H^r$.
The last statement is the motivation of the construction.
Indeed, by Proposition~\ref{junction-exponent}(iii), the optimal trajectory satisfies its hypothesis when $\varphi$ is of class $C^k$ up to the junctions (which is the case when the free arcs are analytic), and then the family generated by the multiple-arc TFC method contains the optimum.
In the linear-quadratic examples of this paper, we have $r = 2$ and $m_i = 3$.
Hence the normal coordinate of every trial trajectory is $C^2$ across the junctions, and the trajectory itself is $C^1$, with a square-integrable second derivative, which is the class in which the linear-quadratic problems of Section~\ref{analysis} are posed.

\begin{remark}[relaxed multiple-arc TFC method]\label{relaxed-remark}
The square in \eqref{F2-free-arc} enforces the nonnegativity of the residual on the free arcs, and tangency alone would not ensure feasibility.
Indeed, if $y(\cdot)$ is parametrized directly, with the contact conditions $y^{(j)}(t_1) = 0$ for $j < p$, then its sign is unrestricted on the free arcs.
For the Bryson-Denham problem, the optimal junction times are not even a local minimizer of this relaxed problem (Lemma~\ref{relaxed-failure}).
Full tangency up to order $m - 1$ restores a local minimum, but the relaxed problem still admits infeasible trajectories with a cost lower than the optimal one (Remark~\ref{relaxed-tangency}).
\end{remark}

\begin{remark}[structure and initialization]\label{structure}
The multiple-arc TFC method takes as given the number of boundary arcs, their order and the junction exponents (as multiple shooting takes the switching structure, see \cite{Trelat2012JOTA}), and it is a local method, which needs an initialization.
In the linear-quadratic class, the two-step procedure used in this paper answers both needs.
In step 1, we run the single-arc TFC method, which requires no structure and whose solution is feasible for every degree, and we locate the boundary arcs where its slack function is small, by a threshold.
This gives the number of boundary arcs and initial junction times, to a few hundredths of the horizon.
In step 2, we run the multiple-arc TFC method from this initialization, and we refine the junction times to the accuracy reported in Section~\ref{examples}.
In all runs of Sections~\ref{bryson} and~\ref{rendezvous}, the initialization provided by step 1 was sufficient for convergence.

The exponent is determined by the junction only up to an even integer.
Indeed, if the true residual is $\delta^{m^\star} \varphi$, then the factor $\delta^m$ with $m \leq m^\star$ and $m^\star - m$ even leaves a factor $\sigma = \delta^{(m^\star - m)/2} \sqrt{\varphi}$, which is as regular as $\sqrt{\varphi}$ and vanishes at the junction, whereas a choice of the wrong parity, or larger than $m^\star$, loses the regularity.
Since $m^\star \geq p$, one of the two candidates $p$ and $p + 1$ is always compatible, and the true exponent is recognized by $\sigma(t_j) \neq 0$, judged by its stability under refinement.

Additional cuts inside a free or a boundary arc, with junction values that impose the continuity of the trajectory and no vanishing factor, do not change the range, and we use them in Section~\ref{reentry} to adapt the polynomial degrees to the arcs.
On the contrary, a cut carrying a vanishing factor forces a contact even where the optimum has none, and a boundary arc of positive duration restricts the family.
Finally, a boundary arc that is allowed to disappear requires a treatment of zero durations, which is not contained in the present formulation.

Of course, step 1 can be replaced by any method providing a plausible structure and initial junction times, for instance a direct method, possibly with an adaptive mesh that locates the junctions \cite{ByczkowskiRao2024}, a continuation on the constraint, as in indirect methods by homotopy \cite{CaillauCotsGergaud2012,HaberkornMartinonGergaud2004}, or a geometric analysis of the extremals \cite{Bonnard2003Reentry,BonnardFaubourgTrelat2005,Trelat2012JOTA}.
For the reentry problem of Section~\ref{reentry}, where the dynamics are not eliminated, step 1 is the single-arc TFC method in a residual form (one constraint is embedded, and the others are imposed at the nodes), started from a trajectory simulated with a constant bank angle, and the junction times are read from the corner of the slack, and not by a threshold.

The procedure does not guarantee that the structure found is the optimal one, and it does not handle a boundary arc that appears or disappears during the iteration.
An outer loop could update the structure: it would detect an interval of persistently small slack inside a free arc, away from the existing junctions, or a boundary arc whose duration becomes very small, and then insert or remove an arc and restart.
These are heuristic indicators, to be checked under refinement, and not certificates of the correct structure.
We have not implemented such an update.
\end{remark}

\section{Exact cost identity and feasible upper bounds}\label{analysis}

\subsection{A linear-quadratic class and its reduced problems}\label{flat}

Throughout this section, we restrict the analysis to a class of flat linear control systems, in which the dynamics are eliminated exactly through the flat output, so that every trial trajectory is admissible.
The state $x = (\xi, \dot \xi, \dots, \xi^{(r-1)})$ is built from an output $\xi(\cdot) : I \to \R^d$, and the control is a linear differential operator with constant coefficients applied to the output, possibly with a known forcing term $f_0(\cdot) \in L^2(I; \R^d)$,
\begin{equation}\label{flat-control}
u = \Lambda \xi - f_0, \qquad \Lambda \xi = \xi^{(r)} + \sum_{j=0}^{r-1} A_j \xi^{(j)}, \qquad A_j \in \R^{d \times d} .
\end{equation}

The cost is $J(\xi) = \frac{1}{2} \int_0^T \vert \Lambda \xi(t) - f_0(t) \vert^2 \dd t$, the endpoint conditions prescribe the state at $t = 0$ and at $t = T$, and the path constraint $a^\top \xi(t) \leq b$ is affine in the output.
For $\xi(\cdot) \in H^r$, the pair $(x, u)$ is admissible for the control system, and the set of admissible trajectories, still denoted by
$$
\adm = \{\xi(\cdot) \in H^r(I; \R^d) \mid L\xi = c, \ a^\top \xi(t) \leq b \text{ on } I\} ,
$$
is a closed convex subset of an affine subspace of $H^r$.

\begin{remark}[flat control systems]
Every controllable linear system is flat, and the class above consists of the systems whose controllability indices are all equal to $r$, with $\xi$ as flat output and the energy of the control as cost \cite{Fliess1995}.
It contains, for instance, integrator chains, linear mechanical systems with as many actuators as degrees of freedom, and the Clohessy-Wiltshire equations of Section~\ref{rendezvous}.
Unequal indices can be handled in the same way, with an operator $\Lambda$ whose rows have different orders, which we do not write out.
For a flat system, linear or not, every trajectory of the flat output is the trajectory of an admissible pair, and this removes the dynamics from the problem.
Linearity makes the cost quadratic in the output and keeps the constraint affine, and these two properties are used in the identity of Proposition~\ref{ritz-theorem}.

A constraint affine in the physical state but involving derivatives of the flat output, such as a velocity bound for a double integrator, requires another construction (Remark~\ref{other-cases}).
Outside this class, the dynamics have to be enforced as residuals, and we examine in Section~\ref{reentry} what remains of the construction in that case.
\end{remark}

The reduced problems are obtained by restricting the free functions to polynomials of degree $N$ on each arc:
\begin{align}
&\min_{g \in \R^{n_1}} J\big(\mathcal{X}[g]\big) , \tag{F1}\label{F1}\\
&\min_{\theta, \ g \in \R^{n_2}} J\big(\mathcal{X}^{(2)}_\theta[g]\big), \qquad 0 < t_1 < t_2 \leq T . \tag{F2}\label{F2}
\end{align}
Both costs are sums of squares, evaluated by a Gauss quadrature that is exact when the integrand is polynomial, and Levenberg-Marquardt iterations can be applied.

We solve \eqref{F1} and then \eqref{F2}, as in the two-step procedure of Remark~\ref{structure}, and the second problem returns a spectrally accurate solution, together with the junction times.
The same division of labor exists between direct and indirect methods \cite{Trelat2012JOTA}, but here both stages live in the same reduced formulation.

Note that two parametrizations whose ranges coincide define the same reduced problem, whatever the free function and the switching functions used to write them.
Indeed, a constrained expression is a parametrization of a set, and only the set matters for the value of \eqref{F1} or \eqref{F2}.
Hence the single-arc TFC method of Definition~\ref{single} and the multiple-arc TFC method differ by their ranges (Propositions~\ref{exact-range} and~\ref{multiple-range}) and by the conditioning of the coordinates that they induce.

\subsection{The identity and its consequences}\label{ritz}

In this section, we establish the following three results: (i) the cost of every admissible trial trajectory is an upper bound of the optimal cost, which follows from admissibility alone; (ii) in the linear-quadratic class, an exact identity decomposes the cost error into a multiplier term and the energy of the discrepancy; (iii) the infimum of the cost over the multiple-arc family converges geometrically when the arcs of the optimal trajectory are analytic.

Let $V = \{h \in H^r \mid Lh = 0\}$.
We assume that $\Lambda$ is coercive on $V$, which holds as soon as the initial state is prescribed, so that $J$ has a unique minimizer $\xi^\star$ on $\adm$, and that the multiplier rule holds at $\xi^\star$: there is a nonnegative Radon measure $\mu$ on $I$, supported on the contact set $\{a^\top \xi^\star = b\}$, such that
\begin{equation}\label{kkt-measure}
\int_0^T (\Lambda \xi^\star - f_0)^\top \Lambda h \dd t + \int_I a^\top h \dd\mu = 0 \qquad \text{for every } h \in V .
\end{equation}

The existence of a Slater point, i.e., of some $\xi_s \in \adm$ with $a^\top \xi_s < b$ on $I$, is sufficient for the multiplier rule \cite{BonnansHermant2009AIHP}, \cite[Section~3.1]{BonnansShapiro2000}.
A Slater point exists for the Bryson-Denham problems, but not for the rendezvous problem, whose target lies on the constraint boundary.
For the rendezvous problem, we read the measure from the shooting solution, and we check \eqref{kkt-measure} directly.
For the problems of Section~\ref{examples}, $\mu$ has a density $\eta$ on the boundary arc and atoms $\nu_j$ at the junctions or at touch points.
For instance, for the Bryson-Denham problem with a boundary arc, we have $\mu = \frac{2}{9\ell^2} (\delta_{3\ell} + \delta_{1 - 3\ell})$, and the density on the arc is zero, whereas the forced variant of Section~\ref{bd-forced} has a nonzero density.

\begin{proposition}[exact error identity]\label{ritz-theorem}
For every $\xi(\cdot) \in \adm$,
\begin{equation}\label{error-identity}
J(\xi) - J(\xi^\star) = \int_I (b - a^\top \xi) \dd\mu + \frac{1}{2} \int_0^T \vert \Lambda(\xi - \xi^\star) \vert^2 \dd t .
\end{equation}
Consequently, for any family $\adm_N \subset \adm$ with $J_N = \inf_{\adm_N} J$: $J_N \geq J^\star$ for every $N$, so that the cost of every trial trajectory is an upper bound of $J^\star$; and if $\adm_N$ contains trajectories $\xi_N$ with $b - a^\top \xi_N = 0$ on $\supp \mu$, then $J_N - J^\star \leq \frac{1}{2} \Vert \Lambda(\xi_N - \xi^\star) \Vert_{L^2}^2$.
\end{proposition}

\begin{proof}
Let $h = \xi - \xi^\star \in V$.
Expanding the quadratic cost, $J(\xi) - J(\xi^\star) = \int (\Lambda \xi^\star - f_0)^\top \Lambda h + \frac{1}{2} \int \vert \Lambda h \vert^2$, and \eqref{kkt-measure} gives $\int (\Lambda \xi^\star - f_0)^\top \Lambda h = -\int a^\top h \dd\mu = \int (b - a^\top \xi) \dd\mu - \int (b - a^\top \xi^\star) \dd\mu$, where the last integral vanishes because $\mu$ is supported on the contact set of $\xi^\star$.
Both terms of \eqref{error-identity} are nonnegative since $\mu \geq 0$ and $\xi \in \adm$.
\end{proof}

This identity is nothing but the quadratic expansion of the cost combined with the multiplier rule.
Its interest is that it separates two contributions: that of staying strictly inside the constraint where the multiplier lives, and the energy of the discrepancy.
For the multiple-arc TFC method with the exact junction times, the first term vanishes identically.

\begin{corollary}[geometric convergence of the multiple-arc TFC method]\label{spectral}
Assume that $\xi^\star$ has the structure of Definition~\ref{multiple-def} with junction times $\theta^\star$ and junction exponents $m_1, m_2 \geq r$, that $\xi^\star$ is real-analytic on a neighborhood of each closed arc, that $y^\star > 0$ on the free arcs except at the junctions, and let $\adm_N$ be the range of $\mathcal{X}^{(2)}_\theta$ with free functions of degree $N$, $\theta$ free, and $\zeta$ encoding the continuity of $z$ up to order $r - 1$.
Then there exist $C > 0$ and $\rho > 1$ such that $J^\star \leq J_N \leq J^\star + C \rho^{-2N}$.
\end{corollary}

\begin{proof}
By Proposition~\ref{multiple-range}, $\adm_N \subset \adm$.
By Proposition~\ref{junction-exponent}(iii) and the positivity of $y^\star$ away from the junctions, $\sigma_1^\star = \sqrt{y^\star / (t_1^\star - t)^{m_1}}$ is analytic and positive on the closed arc $[0, t_1^\star]$, and so are $\sigma_2^\star$ and the components of $z^\star$ on each arc.
Their Chebyshev truncations of degree $N$ converge geometrically in every $C^j$ norm \cite[Chapter~8]{Trefethen2019ATAP}, and the corresponding constrained expressions with the exact junction parameters $\theta^\star$ define $\xi_N \in \adm_N$ with $\Vert \xi_N - \xi^\star \Vert_{H^r} \leq C \rho^{-N}$, since the map $(\sigma, z) \mapsto \xi$ is polynomial in $\sigma$ and linear in $z$ and the retraction onto the transferred data is a bounded linear correction.
Proposition~\ref{ritz-theorem} gives $J(\xi_N) - J^\star \leq C \rho^{-2N}$, and $J_N \leq J(\xi_N)$ because $\theta$ is free.
\end{proof}

Note that the corollary concerns the infimum $J_N$ over the ansatz.
A local solver returns the cost $\widehat J_N$ of a trial trajectory, which is an upper bound of $J_N$ and of $J^\star$.
The convergence curves of Section~\ref{examples} are observations on these computed values, whose last digits are affected by the quadrature and by rounding, and we do not prove the convergence of the iterations.

Without analyticity, the existence of trial trajectories at the exact junctions with $\Vert \xi_N - \xi^\star \Vert_{H^r} \leq \varepsilon_N$ implies $J_N - J^\star \leq C \varepsilon_N^2$, and $C^\infty$ factors lead to a superalgebraic rate, but not to a geometric one in general.

For the single-arc TFC method, the first term of \eqref{error-identity} does not vanish.
Indeed, $s^\star = 0$ on the boundary arc, and thus a polynomial slack $s_N$ satisfying the transferred data with $\Vert s_N - s^\star \Vert_\infty = e_N$ gives $b - a^\top \xi_N = s_N^2 \leq e_N^2$ there, and a first term at most equal to $\mu(I) e_N^2$, which is of the order of $N^{-m}$ by Corollary~\ref{jackson} if $e_N$ follows the best approximation error.
This is an upper bound of one term, and not an estimate of the cost error, whose second term involves the derivatives of $s_N$.
In the experiments, the cost error decays like $N^{-3}$ on the two problems with $m = 3$, but we have no proof of this rate.

\begin{remark}[dual recovery]\label{dual}
Let us consider the second-order problems of this class ($r = 2$).
The Pontryagin maximum principle for state-constrained problems \cite{Clarke2013,HartlSethiVickson1995,Maurer1977,Vinter2000}, in the direct adjoining form, with the Hamiltonian $H = \frac{1}{2} \vert u \vert^2 + p_\xi^\top \dot \xi + p_v^\top (-A_1 \dot \xi - A_0 \xi + u + f_0)$ minimized in $u$, states that the costate $p = (p_\xi, p_v)$ satisfies $p_v = -u^\star$ and $p_\xi = \dot u^\star - A_1^\top u^\star$ on the free arcs, that it jumps by $-\nu_j \nabla c$ at the junctions, with $\nu_j \geq 0$, and that $\dd\mu = \eta \dd t$ with $\eta \geq 0$ on the boundary arcs.
Therefore, on a free arc, the costate is a function of the control and of its derivatives, and the solution of \eqref{F2} directly provides $p$ there.
The jump $\nu$ and the density $\eta$ then follow from the junction conditions and from the invariant $c^{(p)}(x, u) = 0$, differentiated along the arc.

This is the information needed to start an indirect shooting method, and we show in Section~\ref{rendezvous} that a Newton iteration started from the costates recovered from the multiple-arc TFC method converges in two steps.
We do not develop here the correspondence between the stationarity conditions of \eqref{F2} and the junction conditions of the maximum principle.
In the experiments, the multiple-arc TFC method, in which the continuity of the control at the junction is not imposed, produces a control whose jump at $t_1$ decreases with $N$, as required by optimality.
\end{remark}

\begin{remark}[minimum time, box constraints and constrained derivatives]\label{other-cases}
For a minimum-time problem, the feasibility argument gives $T_N \geq T^\star$ for every trial trajectory, but the identity \eqref{error-identity} has no analogue.
When the control is bang-bang, the multiple-arc TFC method reduces to the classical optimization of the switching times.

A bilateral bound $\vert a^\top x(t) \vert \leq b$ cannot be embedded by two independent squared slacks, since the two residuals sum to $2b$.
Bounded constraint sets admit other smooth surjective parametrizations, such as $b \sin \theta$ for an interval, but they lose the single right-inverse correction.

When the constrained quantity is a derivative of the flat output (for instance, for a double integrator with a velocity bound $v \leq v_{\max}$), the endpoint condition on the position becomes the integral constraint $\int_0^T (v_{\max} - s(t)^2) \dd t = x(T) - x(0)$, which is quadratic in the slack and must be kept as an explicit scalar constraint of the reduced problem.

When the dynamics cannot be eliminated, the trial trajectories are not admissible, and neither the upper-bound property nor the identity holds.
Only the exactness of the endpoint and path constraints remains, and we examine this situation in Section~\ref{reentry}.
\end{remark}

The guarantees obtained for the two classes of problems considered hereafter are summarized in Table~\ref{table-guarantees}.
\begin{table}[ht]
\centering\small
\caption{Guarantees for continuous TFC trial trajectories, before numerical quadrature and rounding.}
\label{table-guarantees}
\begin{tabular}{>{\raggedright\arraybackslash}p{4.0cm}>{\raggedright\arraybackslash}p{4.0cm}>{\raggedright\arraybackslash}p{6.0cm}}\toprule
 & linear-quadratic class & reentry, dynamics as residuals\\\midrule
endpoint and path constraints & exact on the whole horizon & endpoint data exact; each path constraint exact only on the arcs carrying its representation\\
dynamical admissibility & guaranteed by elimination & not guaranteed (dynamics imposed as residuals)\\
cost upper bound & yes, for the exact continuous cost & not guaranteed\\
\bottomrule\end{tabular}
\end{table}

\section{Numerical examples}\label{examples}

In this section, we compare the single-arc and multiple-arc TFC methods with direct collocation and indirect shooting on the Bryson-Denham problem and a forced variant of it, on a rendezvous problem with an approach corridor and on the atmospheric reentry of a space shuttle with three state constraints.
We first describe the methods as implemented, then we present the three comparisons, and we end with some remarks.

\paragraph{The TFC methods.}
On each arc, the free functions are expanded in Chebyshev polynomials of degree $N$, and the support functions are the first $q$ Chebyshev polynomials when the arc carries $q$ conditions, so that an arc has $N + 1 - q$ free coefficients.
The cost is evaluated by Gauss-Legendre quadrature, which is exact for polynomial integrands, and the reduced problems \eqref{F1} and \eqref{F2} are solved as nonlinear least squares by the Levenberg-Marquardt method of SciPy \cite{Virtanen2020SciPy}, with analytic Jacobians with respect to the coefficients.
The junction parameters of \eqref{F2} are handled by an outer derivative-free minimization when they reduce to the junction times, and otherwise by finite differences inside the same iteration.

For the single-arc TFC method, the degrees are increased in sequence, each degree being initialized from the previous solution.
The structure and the initial junction times of the multiple-arc TFC method are read from the single-arc solution by a threshold on the slack function, which is step 1 of the procedure of Remark~\ref{structure}.
Neither of the two TFC methods is guaranteed to converge globally, and we check afterwards that the factors $\sigma$ have no zero inside the free arcs.

\paragraph{Direct collocation.}
Our baseline transcriptions use a uniform grid with $N$ intervals and the trapezoidal or the Hermite-Simpson rule in separated form, the path constraint being imposed at the nodes and, for Hermite-Simpson, at the midpoints \cite{Betts2010,HargravesParis1987,Kelly2017}.
The nonlinear programs are built with \texttt{CasADi} \cite{Andersson2019CasADi}, which provides their exact derivatives by automatic differentiation, and they are solved by \texttt{IpOpt} \cite{WachterBiegler2006} with tolerance $10^{-10}$.
A modeling language such as AMPL \cite{FourerGayKernighan2003}, which also provides exact derivatives, could be used as well.
A discrete solution defines a continuous trajectory through its natural interpolant (piecewise linear for the trapezoidal rule, cubic Hermite for Hermite-Simpson), and the reported errors are those of the interpolants, sampled on a fine grid containing the junction times.

\paragraph{Indirect shooting.}
The indirect method takes the arc structure as given and solves the junction conditions of the maximum principle by single shooting, with the routine \texttt{hybrd} of MINPACK (Powell's hybrid method) called through \texttt{fsolve} of SciPy.
Its basin of attraction is measured by the fraction of random initial guesses, drawn uniformly in a box around the solution, from which it converges to the known solution (Sections~\ref{bryson} and~\ref{rendezvous}).

\paragraph{Conventions.}
The tables hereafter report the signed cost error $\widehat J_N - J^\star$ (whose sign distinguishes values above the reference from values below it), the uniform errors on the output and on the control, the junction times when the method produces them, and the number of unknowns of the program actually solved.
The reference solutions are in closed form for the Bryson-Denham problem, computed in extended precision for its forced variant and for the rendezvous problem, and obtained by structured shooting for the reentry problem.

\subsection{The Bryson-Denham problem}\label{bryson}

\paragraph{Three regimes.}
The problem of Bryson, Denham and Dreyfus \cite{BrysonDenhamDreyfus1963}, in the form of \cite[Section~3.11]{BrysonHo1975}, is
\begin{equation}\label{bd}
\min \ \frac{1}{2} \int_0^1 u(t)^2 \dd t, \quad \dot x(t) = v(t), \ \ \dot v(t) = u(t), \quad x(0) = x(1) = 0, \quad v(0) = 1, \ \ v(1) = -1, \quad x(t) \leq \ell ,
\end{equation}
with $\ell > 0$.
Here, the flat output is $\xi = x$, we have $\Lambda = \dd^2 / \dd t^2$, the constraint is of order $p = 2$, and the solution is symmetric with respect to $t = 1/2$, with three regimes.
For $\ell \geq 1/4$, the constraint is inactive, $x^\star(t) = t - t^2$ and $J^\star = 2$.
For $1/6 < \ell < 1/4$, the trajectory touches the bound at $t = 1/2$ with a quadratic contact.
On $[0, 1/2]$, we have $x^\star(t) = t + at^2 + bt^3$, with $a = 12\ell - 4$ and $b = 4 - 16\ell$, the derivative of the control jumps at the touch point, and $J^\star = 96\ell^2 - 48\ell + 8$.
For $\ell < 1/6$, there is a boundary arc $[3\ell, 1 - 3\ell]$.
On $[0, 3\ell]$, we have $x^\star(t) = \ell (1 - (1 - t/3\ell)^3)$ and $u^\star(t) = -\frac{2}{3\ell} (1 - \frac{t}{3\ell})$, whereas $u^\star = 0$ on the arc, and $J^\star = 4/(9\ell)$.
On the free arc, $y^\star(t) = \ell (1 - t/3\ell)^3$, and thus the junction exponent is $m = 3$, as predicted by Proposition~\ref{junction-exponent}.
Recall that, at the touch point, $y^\star / (t - 1/2)^2$ is Lipschitz and not $C^1$ (Section~\ref{junction}).

\paragraph{Single-arc TFC method.}
In the adapted coordinates $y = \ell - x$ (there is no $z$), the slack carries the four transferred conditions $s(0) = \sqrt{\ell}$, $\dot s(0) = -1/(2\sqrt{\ell})$, $s(1) = \epsilon_1 \sqrt{\ell}$, $\dot s(1) = \epsilon_1 / (2\sqrt{\ell})$, with the parity $\epsilon_1 = \pm 1$ of Proposition~\ref{exact-range}.
The control is $u = -2(\dot s^2 + s \ddot s)$, and the reduced problem has $N - 3$ unknowns.
Table~\ref{table-bd} reports the convergence for $\ell = 0.1$ with $\epsilon_1 = 1$.
The cost error decays like $N^{-3}$, as expected from $m = 3$, the position error decays faster, and the control error like $N^{-1}$.
Every computed trajectory satisfies $x \leq \ell$, and every computed cost is above $J^\star$.
The iteration slows down for $N \geq 64$, which is the degeneracy described in Remark~\ref{single-obstructions}(iv).
The threshold applied to the solution with $N = 24$ locates the boundary arc at $[0.2715, 0.7285]$, instead of $[0.3, 0.7]$.

\begin{figure}[h]
\centering
\includegraphics[width=\textwidth]{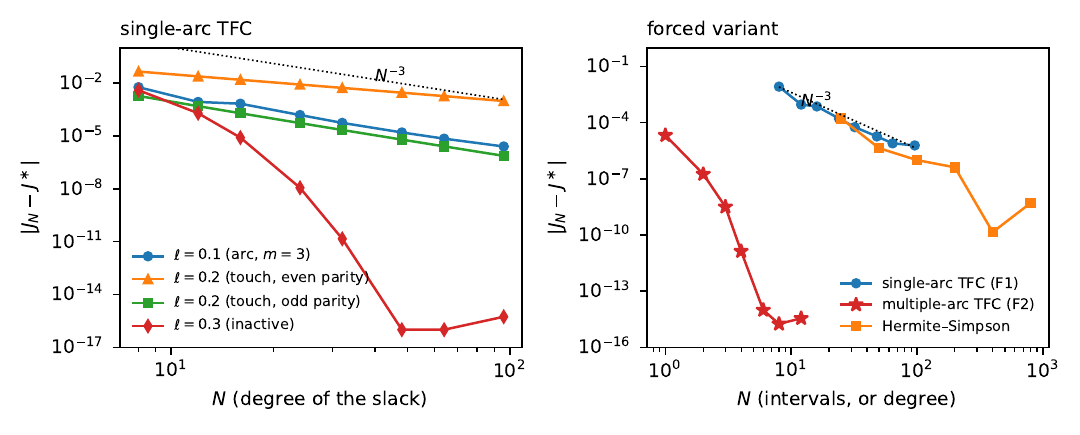}
\caption{Left: absolute cost error of the single-arc TFC method in the three regimes of the Bryson-Denham problem, and for both parities at the touch point.
The dotted line is $N^{-3}$, and a floor of $10^{-16}$ is applied to the values at rounding level.
Right: forced variant, absolute cost error of the two TFC methods and of Hermite-Simpson collocation, against the extended-precision reference, as a function of the degree or of the number of intervals.}
\label{fig-bd-conv}
\end{figure}

The parity matters in the other regimes (Figure~\ref{fig-bd-conv}, left).
In the touch-point regime ($\ell = 0.2$), the odd parity $\epsilon_1 = -1$ lets the slack change sign once, the closure of the ansatz in $H^2$ contains the signed root $(t - 1/2) \sqrt{y^\star / (t - 1/2)^2}$, which is $C^{1,1}$, and the cost error decays like $N^{-3}$.
With the even parity, the slack can only form a cusp there, and the error decays like $N^{-1.5}$.
In the inactive regime ($\ell = 0.3$), the even parity converges geometrically, down to the rounding level for $N \geq 48$, whereas the odd parity cannot converge at all, since it forces a contact.
For instance, at $\ell = 0.26$, its iteration converges to $2.0096$, while $J^\star = 2$, and no degree can reduce this gap (Remark~\ref{parity-gap}).

\begin{remark}[cost of the wrong parity in the inactive regime]\label{parity-gap}
For $\ell > 1/4$, an odd-parity slack $s(\cdot) \in H^2(0, 1)$ vanishes at some $\theta$, so that $x_s = \ell - s^2$ touches the bound there.
Writing $x_s = x_0 + z$, with $x_0(t) = t(1 - t)$ the unconstrained minimizer and $z \in H^2_0(0, 1)$, the cost is $2 + \frac{1}{2} \Vert \ddot z \Vert_{L^2}^2$, and the Cauchy-Schwarz inequality applied to $z(\theta) = \ell - \theta(1 - \theta)$ gives the infimum $2 + 96 (\ell - 1/4)^2$.
This infimum is attained by the touch-point cubic of Section~\ref{bryson}, continued to $\ell > 1/4$, whose signed root has the odd parity, and it is the same over polynomials of any degree.
\end{remark}

\paragraph{Multiple-arc TFC method.}
With $m = 3$, the free arc $[0, t_1]$ carries $y(t) = (t_1 - t)^3 \sigma_1(t)^2$, with the two transferred conditions $\sigma_1(0) = \sqrt{\ell / t_1^3}$ and $\dot \sigma_1(0) = (3\ell / t_1 - 1) / (2 t_1^3 \sigma_1(0))$.
The arc $[t_2, 1]$ is symmetric, and the boundary arc is $x = \ell$.
The exact free arc corresponds to a constant $\sigma_1$, and thus the family contains the optimal trajectory for every $N \geq 1$.
For $N = 1$, the only unknowns are the two junction times, and we obtain $t_1$ within $10^{-8}$ of $3\ell$ (limited by the tolerance of the outer minimization) and the cost within $10^{-14}$ of $J^\star$ (Table~\ref{table-bd}).
This exactness is specific to the example, whose free arcs are cubic, and it says nothing about the rate of convergence when the optimal trajectory does not belong to the family.
This rate is measured on the forced variant below.

In the touch-point regime, treated with two arcs $y = (t_1 - t)^2 \sigma_1^2$ and $y = (t - t_1)^2 \sigma_2^2$ glued by a common junction value, $\sigma^\star$ is not constant, and the cost error decays from $4 \times 10^{-7}$ at $N = 3$ to the rounding level at $N = 12$.

\paragraph{The relaxed method.}
Let us now check on this problem that the contact conditions alone cannot replace the square on the free arcs, as announced in Remark~\ref{relaxed-remark}.
For a candidate entry time $a$, let $x_a(\cdot)$ minimize $\frac{1}{2} \int_0^a \ddot x(t)^2 \dd t$ subject to $x(0) = 0$, $\dot x(0) = 1$ and to the contact conditions $x(a) = \ell$, $\dot x(a) = 0$, without the inequality.
Whatever the degree of the polynomial ansatz, the minimizer is the cubic Hermite interpolant, and we set $V(a) = \frac{1}{2} \int_0^a \ddot x_a(t)^2 \dd t$.
Explicitly, we have
$$
V(a) = \frac{2(a^2 - 3a \ell + 3\ell^2)}{a^3}, \qquad \ddot x_a(a) = \frac{2(a - 3\ell)}{a^2} ,
$$
and, the arc $[t_2, 1]$ being symmetric, the relaxed cost of the free-boundary-free structure with junction times $t_1 \leq t_2$ is $V(t_1) + V(1 - t_2)$.

\begin{lemma}[monotonicity of the relaxed free-arc cost]\label{relaxed-failure}
$V$ is differentiable on $(0, +\infty)$ and $V'(a) = -\frac{1}{2} \ddot x_a(a)^2 \leq 0$.
Hence the relaxed cost is nonincreasing in the length of each free arc; the optimal junction times $t_1^\star = 3\ell$ and $t_2^\star = 1 - 3\ell$, at which $u^\star$ vanishes, are stationary but not locally minimizing for the relaxed problem, and for an entry time slightly beyond $3\ell$ the relaxed cost is smaller and the free arc violates the constraint.
\end{lemma}

\begin{proof}
Write $x = x_a$ and $\partial x = \partial x_a / \partial a$.
Differentiating $x_a(a) = \ell$ and $\dot x_a(a) = 0$ with respect to $a$ gives $\partial x(a) = -\dot x(a) = 0$ and $\partial \dot x(a) = -\ddot x(a)$, while $\partial x(0) = \partial \dot x(0) = 0$.
Then, integrating by parts twice and using $x^{(4)} = 0$,
$$
V'(a) = \frac{1}{2} \ddot x(a)^2 + \int_0^a \ddot x \, \partial \ddot x \dd t = \frac{1}{2} \ddot x(a)^2 + \ddot x(a) \, \partial \dot x(a) = -\frac{1}{2} \ddot x(a)^2 .
$$
At $a^\star = 3\ell$ the interpolant is the exact free arc of the first paragraph and $\ddot x_{a^\star}(a^\star) = u^\star(3\ell) = 0$, as Proposition~\ref{junction-exponent} predicts with $m = 3$; $V'$ is negative on both sides of $a^\star$, so $V$ keeps decreasing as $a$ increases.
For $a > a^\star$, $\ddot x_a(a) > 0$ and $x_a(a - h) = \ell + \frac{1}{2} \ddot x_a(a) h^2 + \mathrm{O}(h^3) > \ell$ for small $h > 0$: the cubic approaches $\ell$ from above.
\end{proof}

Minimized over the junction times with $t_1 \leq t_2$, the relaxed method behaves as predicted by the lemma: the boundary arc collapses, and for $\ell = 0.1$ the minimizer is $t_1 = t_2 = 0.25$, with $J = 3.98$, below $J^\star = 4.44$, and its second free arc exceeds the bound by $0.045$ (Table~\ref{table-bd}).

\begin{remark}[keeping the full tangency]\label{relaxed-tangency}
In the relaxed method, we have imposed only the contact conditions.
Imposing the full tangency of order $m - 1 = 2$, i.e., $\ddot x(a) = 0$ in addition to the contact conditions, on polynomials of degree $D \geq 4$, we get the relaxed cost $W_D(a) = V(a) + 2(a - 3\ell)^2 / \big(a^3 (D - 3)(D + 1)\big)$, with $W_D''(a^\star) > 0$.
Hence the optimal junction times become a strict local minimizer, with a curvature that vanishes as $D \to +\infty$.
But the global infimum remains below $J^\star = 4.444$ and is infeasible ($4.200$ for $D = 4$ and $4.040$ for $D = 6$, with a collapsed boundary arc), so that neither form of the relaxation is safe.
\end{remark}

The lemma is the primal counterpart of the sign conditions of the Pontryagin maximum principle.
In the indirect method, the multiplier density and the junction jumps must be nonnegative, and similarly, in the reduced primal problem, the feasibility of the free arcs must be enforced, since otherwise the minimization over the junction times shortens the boundary arc.

\paragraph{Collocation and shooting.}
The junctions $0.3$ and $0.7$ are nodes of every grid whose $N$ is a multiple of ten.
On such grids, the Hermite-Simpson transcription represents the piecewise-cubic solution and reaches the solver tolerance.
On the other grids, its cost error decays like $N^{-3}$ and stagnates at $3 \times 10^{-8}$, and the trapezoidal transcription converges like $N^{-2}$ from above.
Therefore, on this benchmark, a collocation whose grid contains the junctions is not beaten on accuracy, and the single-arc TFC method has the largest control errors of all methods.
The TFC methods provide in addition the exact feasibility of the continuous trajectory, the upper bound and, for the multiple-arc TFC method, the junction time.

The indirect method has six unknowns ($p_x(0)$, $p_v(0)$, $t_1$, $t_2$ and the entry jumps $\nu_1$, $\nu_2$ of the costate) and six junction and terminal equations, solved in closed form since the arcs are explicit.
Its basin of attraction is wide: $377$ of $500$ random initial guesses converge to the solution.

\begin{table}[h]
\centering\small
\caption{Bryson-Denham problem, $\ell = 0.1$: cost error, junction-time error ($t_1^\star = 0.3$), uniform errors on the position and on the control, and number of unknowns, for the indirect shooting method (exact, six unknowns), the single-arc TFC method (degree $N$), the multiple-arc TFC method (degree $N$ per arc), the relaxed method ($t_1 = t_2 = 0.25$, boundary arc collapsed) and the trapezoidal and Hermite-Simpson collocations ($N$ intervals, bound at nodes and midpoints, errors of the interpolants, junctions off the grid unless stated).}
\label{table-bd}
\setlength{\tabcolsep}{4pt}
\begin{tabular}{lrrrrrr}\toprule
method & $N$ & $\widehat J_N - J^\star$ & $t_{1,N} - t_1^\star$ & $\Vert x_N - x^\star \Vert_\infty$ & $\Vert u_N - u^\star \Vert_\infty$ & unknowns\\
\midrule
indirect shooting & - & $0$ & $0$ & $0$ & $0$ & 6\\
\midrule
single-arc TFC & 8 & $5.8 \times 10^{-3}$ & - & $9.0 \times 10^{-4}$ & $3.4 \times 10^{-1}$ & 5\\
 & 16 & $6.7 \times 10^{-4}$ & - & $4.9 \times 10^{-5}$ & $1.6 \times 10^{-1}$ & 13\\
 & 32 & $5.5 \times 10^{-5}$ & - & $4.4 \times 10^{-6}$ & $7.1 \times 10^{-2}$ & 29\\
 & 64 & $6.9 \times 10^{-6}$ & - & $4.9 \times 10^{-7}$ & $3.5 \times 10^{-2}$ & 61\\
 & 96 & $2.5 \times 10^{-6}$ & - & $2.0 \times 10^{-7}$ & $2.6 \times 10^{-2}$ & 93\\
\midrule
multiple-arc TFC & 1 & $-8.9 \times 10^{-16}$ & $-5.5 \times 10^{-9}$ & $<10^{-8}$ & - & 2\\
 & 6 & $-6.2 \times 10^{-15}$ & $-4.0 \times 10^{-9}$ & $<10^{-8}$ & - & 12\\
\midrule
relaxed method & 3 & $-0.46$ & $-0.05$ & $0.045$ (overshoot) & - & 2\\
\midrule
trapezoidal collocation & 85 & $6.6 \times 10^{-3}$ & - & $1.1 \times 10^{-4}$ & $1.1 \times 10^{-1}$ & 258\\
 & 325 & $4.6 \times 10^{-4}$ & - & $7.8 \times 10^{-6}$ & $3.3 \times 10^{-2}$ & 978\\
 & 1285 & $3.0 \times 10^{-5}$ & - & $2.2 \times 10^{-6}$ & $8.6 \times 10^{-3}$ & 3858\\
\midrule
Hermite-Simpson collocation & 45 & $4.2 \times 10^{-5}$ & - & $2.7 \times 10^{-6}$ & $6.2 \times 10^{-2}$ & 183\\
 & 165 & $8.6 \times 10^{-7}$ & - & $1.3 \times 10^{-6}$ & $1.7 \times 10^{-2}$ & 663\\
 & 645 & $2.6 \times 10^{-8}$ & - & $2.8 \times 10^{-6}$ & $4.9 \times 10^{-3}$ & 2583\\
\quad junctions on the grid & 50 & $1.8 \times 10^{-9}$ & - & $1.4 \times 10^{-6}$ & $1.1 \times 10^{-4}$ & 203\\
\bottomrule\end{tabular}
\end{table}

\paragraph{A forced variant.}\label{bd-forced}
In order to remove the polynomial exactness specific to the example, while keeping a closed-form reference, we add a forcing term to the dynamics,
\begin{equation}\label{bd-forced-eq}
\dot v(t) = u(t) + \varepsilon \cos(\omega t), \qquad \varepsilon = \frac{1}{2}, \quad \omega = 2\pi ,
\end{equation}
with the data of \eqref{bd} and $\ell = 0.1$ (in the notation of \eqref{flat-control}, $f_0(t) = \varepsilon \cos(\omega t)$).
The structure is unchanged (free-boundary-free, with $m = 3$, and symmetric with respect to $t = 1/2$), but the free arcs are no longer polynomial, and the boundary control $u = -\varepsilon \cos(\omega t)$ is not zero.
On the free arc $[0, a]$, the control is affine, $u = A + Bt$, the continuity of the control and the tangency at $a$ are linear in $(A, B)$, and $x(a) = \ell$ is a scalar equation for $a$.
Its root, computed in extended precision, is $a = 0.302657815087208$, with $J^\star = 5.24336142728517$.
The multiplier measure of Proposition~\ref{ritz-theorem} is
$$
\mu = \nu (\delta_a + \delta_{1-a}) - \varepsilon \omega^2 \cos(\omega t) \mathbf{1}_{[a, 1-a]} \dd t
$$
with $\nu = B - \varepsilon \omega \sin(\omega a) = 21.58$.
The density on the arc is positive, and thus both terms of the identity \eqref{error-identity} are active.

Table~\ref{table-bd-forced} and Figure~\ref{fig-bd-conv} report the results.
The multiple-arc TFC method is no longer exact, and it converges geometrically to the rounding level with twelve unknowns, the cost errors being positive down to that level, and the junction time is obtained to $3 \times 10^{-8}$.
The single-arc TFC method follows the $N^{-3}$ trend up to $N = 64$, and then stagnates.
Hermite-Simpson collocation converges irregularly and with either sign, since the position of the junctions in the grid varies with $N$, and it violates the bound slightly between the nodes.
At $N = 800$, its control error is $2.4 \times 10^{-3}$, against $6.5 \times 10^{-7}$ for the multiple-arc TFC method with twelve unknowns.

\begin{table}[h]
\centering\small
\caption{Forced Bryson-Denham problem \eqref{bd-forced-eq}, $\ell = 0.1$: cost error with respect to the extended-precision reference, junction-time error, uniform errors on the position and the control, and number of unknowns, for the single-arc TFC method (degree $N$), the multiple-arc TFC method (degree $N$ per arc) and Hermite-Simpson collocation ($N$ intervals).
The errors of the TFC methods are measured on a grid containing the junction times, and those of the collocation are the errors of its interpolants.}
\label{table-bd-forced}
\setlength{\tabcolsep}{4pt}
\begin{tabular}{lrrrrrr}\toprule
method & $N$ & $\widehat J_N - J^\star$ & $t_{1,N} - a$ & $\Vert x_N - x^\star \Vert_\infty$ & $\Vert u_N - u^\star \Vert_\infty$ & unknowns\\
\midrule
single-arc TFC & 8 & $8.4 \times 10^{-3}$ & - & $5.1 \times 10^{-4}$ & $3.2 \times 10^{-1}$ & 5\\
 & 16 & $7.4 \times 10^{-4}$ & - & $4.8 \times 10^{-5}$ & $1.6 \times 10^{-1}$ & 13\\
 & 32 & $5.7 \times 10^{-5}$ & - & $4.9 \times 10^{-6}$ & $7.2 \times 10^{-2}$ & 29\\
 & 64 & $7.8 \times 10^{-6}$ & - & $8.3 \times 10^{-7}$ & $3.5 \times 10^{-2}$ & 61\\
 & 96 & $6.1 \times 10^{-6}$ & - & $5.4 \times 10^{-7}$ & $3.4 \times 10^{-2}$ & 93\\
\midrule
multiple-arc TFC & 1 & $2.1 \times 10^{-5}$ & $1.1 \times 10^{-3}$ & $1.8 \times 10^{-5}$ & $2.3 \times 10^{-2}$ & 2\\
 & 2 & $1.8 \times 10^{-7}$ & $8.9 \times 10^{-5}$ & $1.1 \times 10^{-6}$ & $2.6 \times 10^{-3}$ & 4\\
 & 3 & $3.2 \times 10^{-9}$ & $-1.7 \times 10^{-5}$ & $9.0 \times 10^{-8}$ & $3.6 \times 10^{-4}$ & 6\\
 & 4 & $1.4 \times 10^{-11}$ & $-1.0 \times 10^{-6}$ & $4.4 \times 10^{-9}$ & $2.6 \times 10^{-5}$ & 8\\
 & 6 & $-8.0 \times 10^{-15}$ & $3.0 \times 10^{-8}$ & $1.2 \times 10^{-10}$ & $6.5 \times 10^{-7}$ & 12\\
\midrule
Hermite-Simpson collocation & 50 & $-4.3 \times 10^{-6}$ & - & $2.2 \times 10^{-6}$ & $3.5 \times 10^{-2}$ & 203\\
 & 200 & $4.1 \times 10^{-7}$ & - & $4.8 \times 10^{-7}$ & $1.5 \times 10^{-2}$ & 803\\
 & 400 & $-1.4 \times 10^{-10}$ & - & $3.7 \times 10^{-8}$ & $2.6 \times 10^{-3}$ & 1603\\
 & 800 & $5.1 \times 10^{-9}$ & - & $1.6 \times 10^{-8}$ & $2.4 \times 10^{-3}$ & 3203\\
\bottomrule\end{tabular}
\end{table}

\subsection{Rendezvous with an affine approach corridor}\label{rendezvous}

The relative motion of a chaser with respect to a target on a circular orbit is described, in the orbital plane and with time in units of the inverse mean motion, by the Clohessy-Wiltshire equations \cite{ClohessyWiltshire1960}, where $x$ is radial and $y$ along-track:
$$
\ddot x(t) = 3x(t) + 2\dot y(t) + u_x(t), \qquad \ddot y(t) = -2\dot x(t) + u_y(t) .
$$
The flat output is the relative position $\xi = (x, y)$, with $\Lambda \xi = \ddot \xi + A_1 \dot \xi + A_0 \xi$, $A_1 = \big(\begin{smallmatrix} 0 & -2 \\ 2 & 0 \end{smallmatrix}\big)$ and $A_0 = \big(\begin{smallmatrix} -3 & 0 \\ 0 & 0 \end{smallmatrix}\big)$.

The chaser starts at rest at $(x_0, y_0) = (0.5, -3)$, and it must reach the origin at rest at $T = 2$ with minimal energy $\frac{1}{2} \int_0^T \vert u(t) \vert^2 \dd t$, inside a cone of half-angle $\alpha = 15^\circ$ around the negative along-track axis (a glide-slope constraint of the kind used in proximity operations, see \cite{AcikmesePloen2007,LuLiu2013}).
Since the unconstrained transfer leaves the corridor only through its lower edge, we keep the single affine constraint $c(x, y) = -x + y \tan \alpha \leq 0$, of order $p = 2$, and we check afterwards that the other edge is not touched before the final point, where the two edges meet.
The optimal trajectory follows the free arc until $t_1^\star = 1.4767448$, then slides along the lower edge down to the apex, so that the boundary arc ends at the final time.
Along the boundary arc, the motion has one degree of freedom (the coordinate $z$ along the edge), whose control is optimized.
As already mentioned, the control is continuous at the junction and the derivative of its normal component is not, so that $m = 3$.
The constraint is active at the final point, which is the case of an active endpoint considered in Section~\ref{adapted}.

The reference solution is computed by shooting in the direct adjoining form, with $X = (x, y, \dot x, \dot y)$ and the costate $p$.
On both arcs, the control is $u = -(p_3, p_4)$ and the dynamics of $(X, p)$ are linear, the costate jumps by $-\nu \nabla c$ at the entry, and the unknowns $p(0)$ and $t_1$ are determined by the contact conditions and the continuity of the control at $t_1$, and by the two terminal conditions along the edge.
We obtain
$$
p(0) = (9.26136078, -7.48830176, 6.79702412, -0.61401373), 
\qquad 
t_1 = 1.47674481,
\qquad 
\nu = 10.269489, 
$$
and the multiplier density increases from $18.7$ at the entry to $43.8$ at the apex, so that the multiplier rule \eqref{kkt-measure} holds with a nonnegative measure.
The cost $J^\star = \frac{1}{2} p(0)^\top X(0) = 13.547792827$ is confirmed by an extended-precision solution of the shooting equations.

The adapted coordinates are $y_c = x - y \tan \alpha = s^2$ and $z = (x \tan \alpha + y) / \sqrt{1 + \tan^2 \alpha}$.
The transferred data are $s(0) = \sqrt{x_0 - y_0 \tan \alpha}$, $\dot s(0) = 0$ and $s(T) = 0$, and the single-arc TFC method has $2N - 5$ unknowns.
The multiple-arc TFC method carries $y_c = (t_1 - t)^3 \sigma^2$ on $[0, t_1]$ and $y_c \equiv 0$ on $[t_1, T]$, with the junction values $z(t_1)$, $\dot z(t_1)$ as unknowns besides $t_1$ and the $3N - 7$ coefficients.
It is initialized from the single-arc solution with $N = 32$.

\begin{figure}[h]
\centering
\includegraphics[width=\textwidth]{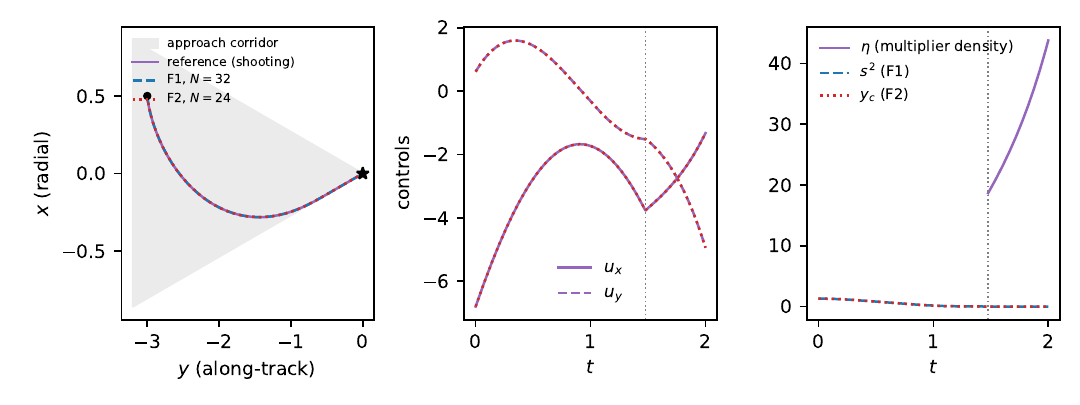}\\[1ex]
\includegraphics[width=0.5\textwidth]{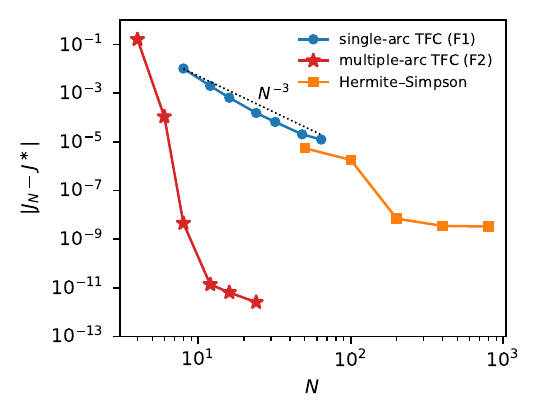}
\caption{Rendezvous problem.
Top left: corridor (shaded), reference trajectory and the two TFC methods in the orbital plane, from the initial point (dot) to the target (star).
Top middle: controls of the reference and of the multiple-arc TFC method, with the entry time (dotted line).
Top right: multiplier density $\eta$ on the boundary arc, squared slack of the single-arc TFC method and constraint residual $y_c$ of the multiple-arc TFC method.
Bottom: absolute cost error of the single-arc TFC method, of the multiple-arc TFC method and of Hermite-Simpson collocation.
The dotted line is $N^{-3}$.}
\label{fig-cw}
\end{figure}

Table~\ref{table-cw} and Figure~\ref{fig-cw} compare the two TFC methods and Hermite-Simpson collocation with the shooting solution.
The single-arc TFC method converges like $N^{-3}$, and its threshold places the entry six percent too early.
The multiple-arc TFC method corrects the entry time, and it converges geometrically on the cost, down to $10^{-11}$ with $32$ unknowns, with positive cost errors and feasible trajectories.
But the junction time stagnates at $5 \times 10^{-6}$ from $N = 8$ on, and so do the position and control errors, because the reduced cost is flat in $t_1$ around its minimum and the iteration stops on its cost tolerance in this valley.
The jump of the control at the junction, which is not imposed, decreases from $1.3$ at $N = 4$ to $10^{-7}$ at $N = 24$, as required by optimality.

Hermite-Simpson collocation reaches a cost error of a few $10^{-9}$ at $N = 800$, with $6406$ unknowns.
With $N = 12$, the multiple-arc TFC method is two orders of magnitude more accurate on the cost and sixty times more accurate on the control, with two hundred times fewer unknowns, and it provides the junction time.
A comparison with an adaptive collocation and a variable junction, motivated by \cite{ByczkowskiRao2024}, is given at the end of this section.

The costate recovered from the solution with $N = 24$ by Remark~\ref{dual} has a shooting residual of $1.4 \times 10^{-3}$, and two Newton steps bring it to $2.4 \times 10^{-14}$.
For comparison, the basin of attraction of the shooting method, measured with initial costates drawn in $[-30, 30]^4$ and random entry times, contains $132$ of $200$ trials.

\begin{table}[h]
\centering\small
\caption{Rendezvous problem: cost error with respect to the shooting solution, entry time, uniform errors on the position and on the control, jump of the control at the junction, and number of unknowns, for the indirect shooting method, the single-arc TFC method (degree $N$), the multiple-arc TFC method (degree $N$ per arc) and Hermite-Simpson collocation ($N$ intervals).
The values $t_1^\star = 1.47674481$ and $J^\star = 13.547792827$ are those of the shooting solution.
For the single-arc TFC method, the entry time is the threshold estimate, and the errors of the collocation are those of its interpolants.}
\label{table-cw}
\setlength{\tabcolsep}{3pt}
\begin{tabular}{lrrrrrrr}\toprule
method & $N$ & $\widehat J_N - J^\star$ & $t_{1,N}$ & $\Vert \xi_N - \xi^\star \Vert_\infty$ & $\Vert u_N - u^\star \Vert_\infty$ & jump & unknowns\\
\midrule
indirect shooting (reference) & - & $0$ & $1.47674481$ & $0$ & $0$ & $0$ & 5\\
\midrule
single-arc TFC & 8 & $1.0 \times 10^{-2}$ & $1.414$ & $2.2 \times 10^{-3}$ & $2.7 \times 10^{-1}$ & - & 11\\
 & 32 & $6.6 \times 10^{-5}$ & $1.387$ & $1.6 \times 10^{-5}$ & $5.9 \times 10^{-2}$ & - & 59\\
 & 64 & $1.2 \times 10^{-5}$ & $1.387$ & $5.3 \times 10^{-6}$ & $3.7 \times 10^{-2}$ & - & 123\\
\midrule
multiple-arc TFC & 4 & $1.6 \times 10^{-1}$ & $1.39304694$ & $3.2 \times 10^{-2}$ & $1.2$ & $1.3$ & 8\\
 & 6 & $1.1 \times 10^{-4}$ & $1.47647067$ & $3.1 \times 10^{-4}$ & $4.0 \times 10^{-2}$ & $4.1 \times 10^{-2}$ & 14\\
 & 8 & $4.5 \times 10^{-9}$ & $1.47675215$ & $8.1 \times 10^{-7}$ & $3.0 \times 10^{-4}$ & $2.7 \times 10^{-4}$ & 20\\
 & 12 & $1.4 \times 10^{-11}$ & $1.47675030$ & $8.8 \times 10^{-8}$ & $5.6 \times 10^{-5}$ & $1.6 \times 10^{-6}$ & 32\\
 & 24 & $2.5 \times 10^{-12}$ & $1.47674936$ & $7.3 \times 10^{-8}$ & $4.4 \times 10^{-5}$ & $1.3 \times 10^{-7}$ & 68\\
\midrule
Hermite-Simpson collocation & 50 & $-5.4 \times 10^{-6}$ & - & $8.2 \times 10^{-5}$ & $2.8 \times 10^{-2}$ & - & 406\\
 & 200 & $6.9 \times 10^{-9}$ & - & $4.7 \times 10^{-6}$ & $1.4 \times 10^{-2}$ & - & 1606\\
 & 800 & $3.3 \times 10^{-9}$ & - & $2.9 \times 10^{-7}$ & $3.5 \times 10^{-3}$ & - & 6406\\
\bottomrule\end{tabular}
\end{table}

\paragraph{Adaptive collocation with a free entry time.}
Let us now compare the multiple-arc TFC method with a two-phase Hermite-Simpson transcription in which the entry time is an unknown.
Both methods start from the same degree-4 constrained curve, built from the single-arc solution with $N = 32$, whose estimated entry time is $t_1 = 1.3865$.
The collocation uses normal and tangential coordinates, and it eliminates the normal coordinate on the boundary arc.
At the entry, it imposes a zero normal position, a zero normal velocity and a zero value of the normal acceleration variable of the scheme, consistently with the junction exponent $m = 3$ used by the multiple-arc TFC method.
Starting from $8 + 4$ intervals, it bisects every interval on which the product of the length by the maximal scaled differential defect at seven interior Gauss points exceeds $\varepsilon$.
The shooting solution is used only to evaluate the errors, and neither for the initialization nor for the refinement.
The results are reported in Table~\ref{table-cw-adaptive}.

\begin{table}[ht]
\centering\small\setlength{\tabcolsep}{5pt}
\caption{Rendezvous problem: multiple-arc TFC method and adaptive Hermite-Simpson collocation, started from a common initialization.
Here, $n$ is the number of unknowns, $\Delta J$ is computed with the continuous energy of the reconstructed control, and $e_\xi$ and $v_c$ are the position error and the corridor violation of the interpolant, sampled at $20001$ points at least, including the junctions.
Times are indicative cumulative seconds on the same machine, including the refinement of degrees or meshes and the construction of the model, and excluding the common initialization (0.08~s) and the validation against the reference.}
\label{table-cw-adaptive}
\begin{tabular}{llrrrrrr}\toprule
method & parameter & $n$ & $\vert \Delta J \vert$ & $\vert \Delta t_1 \vert$ & $e_\xi$ & $v_c$ & time\\\midrule
multiple-arc TFC & $N = 8$ & 20 & $4.5 \times 10^{-9}$ & $7.7 \times 10^{-6}$ & $8.3 \times 10^{-7}$ & $<10^{-15}$ & 5.30\\
 & $N = 12$ & 32 & $1.9 \times 10^{-11}$ & $6.8 \times 10^{-6}$ & $7.2 \times 10^{-8}$ & $<10^{-15}$ & 5.41\\
\midrule
adaptive Hermite-Simpson & $\varepsilon = 10^{-8}$ & 938 & $1.2 \times 10^{-7}$ & $2.3 \times 10^{-5}$ & $6.5 \times 10^{-8}$ & $<10^{-15}$ & 11.29\\
 & $\varepsilon = 10^{-10}$ & 2314 & $6.7 \times 10^{-10}$ & $3.8 \times 10^{-6}$ & $2.5 \times 10^{-10}$ & $<10^{-15}$ & 18.73\\
\bottomrule\end{tabular}
\end{table}

The final collocation mesh has $256 + 64$ intervals.
It provides a more accurate position and a more accurate entry time than the multiple-arc TFC runs of the table, whereas the multiple-arc TFC method reaches a smaller cost error with fewer unknowns.
Note that the small corridor violation of the collocation interpolant does not imply exact dynamical admissibility: integrating independently its reconstructed control, we get a maximal sampled violation of $1.8 \times 10^{-11}$ and a terminal-state residual of $3.6 \times 10^{-11}$.
This comparison concerns only the refinement for a common arc structure: it does not address the detection of the structure, and it does not show that the TFC method is superior to adaptive methods in general.

\subsection{Atmospheric reentry with three state constraints}\label{reentry}

The last example is outside the class of Section~\ref{flat} in every respect: the dynamics are nonlinear, the control is an angle, the final time is free and there are three nonlinear state constraints.
We keep the exact representation of each constraint on the arcs adjacent to its boundary arcs, by a nonlinear pullback instead of the affine change of coordinates \eqref{adapted-coordinates}, the multiple-arc TFC method with the junction times as unknowns, and Proposition~\ref{junction-exponent}.
But the dynamics are not eliminated anymore, and we lose the admissibility of the trial trajectories, the identity \eqref{error-identity} and the bounds.
The single-arc TFC method can still be used in a residual form, with one constraint embedded and the dynamics enforced at the nodes, and it still serves as step 1.
In the computations reported below, step 1 estimates the arc structure and rough junction times, without any prescribed structure of active arcs and without any initial costate, and step 2 refines the solution for the structure so obtained.

\paragraph{Problem.}
We consider the reentry of a space shuttle studied in \cite{Bonnard2003Reentry,BonnardFaubourgTrelat2005,BonnardTrelat2002COCV,Trelat2003DCDS}, with the state $(r, v, \gamma, L, l, \chi)$ (radius, speed, flight-path angle, latitude, longitude and azimuth) and the bank angle $\mu$ as control.
With $g(r) = g_0 / r^2$, $\rho(r) = \rho_0 \exp(-(r - r_T)/h_s)$, $k_D = SC_D/2m$, $k_L = SC_L/2m$ and the Earth rotation rate $\Omega$, the equations of motion are
\begin{equation}\label{reentry6-dyn}
\begin{aligned}
\dot r &= v \sin \gamma, \\
\dot v &= -g \sin \gamma - k_D \rho v^2 + \Omega^2 r \cos L (\sin \gamma \cos L - \cos \gamma \sin L \cos \chi), \\
\dot \gamma &= \cos \gamma \Big(\frac{v}{r} - \frac{g}{v}\Big) + k_L \rho v \cos \mu + 2\Omega \cos L \sin \chi + \Omega^2 \frac{r}{v} \cos L (\cos \gamma \cos L + \sin \gamma \sin L \cos \chi), \\
\dot L &= \frac{v}{r} \cos \gamma \cos \chi, \qquad \dot l = \frac{v}{r} \frac{\cos \gamma \sin \chi}{\cos L}, \\
\dot \chi &= k_L \rho v \frac{\sin \mu}{\cos \gamma} + \frac{v}{r} \cos \gamma \tan L \sin \chi + 2\Omega (\sin L - \tan \gamma \cos L \cos \chi) + \Omega^2 \frac{r}{v} \frac{\sin L \cos L \sin \chi}{\cos \gamma} ,
\end{aligned}
\end{equation}
where the dependence on $t$ is omitted.

The data, the boundary conditions and the incidence law are those of \cite[Table~1 and Appendix~A]{BonnardFaubourgTrelat2005}, with the initial longitude free.
The initial altitude is $119.82$~km and the final altitude is $15$~km, with $l(0)$, $\chi(0)$, $\gamma(t_f)$, $\chi(t_f)$ and $t_f$ free.
The longitude, which is a cyclic variable, is dropped from the state and recovered afterwards.
The drag and lift coefficients are tabulated as functions of the Mach number and of the incidence (see \cite{BonnardFaubourgTrelat2005}).
Composed with the incidence law, they become functions of the Mach number with slope discontinuities at the table breakpoints and at Mach $2$ and $10$.
We smooth them by a Gaussian of standard deviation $0.15$ in Mach, the deviation from the tables being below $5 \times 10^{-3}$, and this smoothed model is the benchmark used by every method below.

The cost is the total thermal flux $\int_0^{t_f} \varphi \dd t$, and the three state constraints are
$$
\varphi = C_q \sqrt{\rho} v^3 \leq \varphi_{\max}, \qquad a_n = \frac{S}{2m} \rho v^2 \sqrt{C_D^2 + C_L^2} \leq a_{\max}, \qquad P = \frac{1}{2} \rho v^2 \leq P_{\max} ,
$$
with $C_q = 1.705 \times 10^{-4}$, $\varphi_{\max} = 717300$~W\,m$^{-2}$, $a_{\max} = 29.34$~m\,s$^{-2}$ and $P_{\max} = 25$~kPa.

\paragraph{Pullbacks, boundary arcs and junction exponents.}
The three constraints are functions of $(r, v)$, and they are of order two.
Indeed, denoting by $c$ any of them, normalized by its bound, $\dot c$ depends on the state and not on $\mu$, and $\ddot c = A + B \cos \mu + C \sin \mu$, with $A, B, C$ functions of the state.
On a boundary arc of $c$, the tangency condition $\dot c = 0$ is solved for the flight-path angle, $\gamma = \Gamma_c(r, v, L, \chi)$, and $\ddot c = 0$ is solved for the bank angle, $\mu = \theta \pm \arccos(-A/R)$, with $R = \sqrt{B^2 + C^2}$ and $\theta = \operatorname{atan2}(C, B)$.
This gives two lateral branches, which differ by the sign of $\sin \mu$.

Each constraint is solved for the speed, $v = V_c(r, y)$, with the residual $y = 1 - c \geq 0$.
On a free arc adjacent to a boundary arc of $c$, the state is represented by $(r, y, \gamma, L, \chi)$, with $y = \delta^m \sigma^2$, where $\delta$ is the distance to the junction, so that $c \leq 1$ holds identically there.
On the boundary arc, the state is represented by $(r, L, \chi)$, the internal degrees of freedom being the lateral motion.
At every junction, the bank angle jumps between $0$ and the boundary value, and the jump of $\ddot c$ does not vanish in the computed solutions.
Hence, by Proposition~\ref{junction-exponent}(i), the junction exponent is $m = p = 2$ at every junction.

Since the dynamics are not eliminated (the system is not known to be flat, and anyway the constraints are not affine in the state), the multiple-arc TFC method is a boundary value solver.
The residuals are the equations \eqref{reentry6-dyn} at the Gauss nodes of each arc, written in the coordinates of the arc, and the unknowns are the junction times, $\chi(0)$, the junction values and the coefficients.
The least-squares problem is solved by a trust-region method with a finite-difference Jacobian, and the cost is evaluated afterwards by quadrature.
The boundary bank angle keeps sharp variations near the breakpoints of the original tables, with a bump of $20^\circ$ on the acceleration arc near the original kink of the incidence law at Mach $10$, and smaller ones elsewhere.
Since a polynomial of moderate degree cannot follow such a control, we cut the arcs at these points, which is allowed at no cost by the multiple-arc formulation (Remark~\ref{structure}).

\paragraph{Step 1: the single-arc TFC method in residual form.}
We use one Chebyshev polynomial of degree $N$ per state component on $[0, 1]$, the endpoint conditions being embedded by constrained expressions.
The flux constraint is embedded by its pullback with a polynomial slack $\sigma$, here solved for the altitude, $r = R_\varphi(v, \sigma)$, with $\varphi / \varphi_{\max} = e^{-\sigma^2}$.
This form agrees with $1 - \sigma^2$ near the boundary, and it keeps the pullback regular where the flux is a few percent of its bound.
The equations \eqref{reentry6-dyn} are imposed as residuals at the $N + 1$ Gauss nodes, with $(\cos \mu, \sin \mu)$ at the nodes as unknowns on the unit circle, and the acceleration and pressure constraints are imposed as inequalities at the nodes.
The cost $\varphi_{\max} t_f \int_0^1 e^{-\sigma^2} \dd \tau$ is minimized by \texttt{IpOpt}.
In other words, the method is a pseudo-spectral collocation with one constraint embedded, and its trial trajectories are not admissible.

The initial guess is the trajectory obtained by integrating \eqref{reentry6-dyn} from the initial state, with a constant bank angle of $60^\circ$ and an initial azimuth of $90^\circ$ (due east), until the final altitude is reached.
This flight lasts $1061$~s, and it satisfies neither the constraints nor the terminal conditions.
The endpoint conditions are then imposed by the constrained expressions.
We take $N = 30$, $40$, $60$ and $80$, each degree being started from the previous solution.
The iteration converges at every degree (in $8$, $9$, $48$ and $251$~s) to a structure in which the three constraints are active in succession (flux, acceleration and pressure), joined by short arcs with $\mu = 0$, after a bang arc with $\mu = \pi$ and an arc with $\mu = 0$.
With $N = 80$ ($561$ unknowns), we obtain $t_f = 877.0$~s and $J = 3.6159 \times 10^8$~J\,m$^{-2}$, and the constraints are exceeded by at most $0.8\%$ between the nodes.

The junction times are read by a hinge fit from the corner of the slack (which, by Proposition~\ref{junction-exponent}, is linear in the distance to the junction on the free side and small on the boundary arc), and from the corner of $\sqrt{1 - c}$ for the two other constraints.
The lateral branch of each boundary arc is given by the majority sign of $\sin \mu$ at its nodes.
On the flux arc, this sign oscillates from node to node, the lateral motion being a nearly degenerate direction.
Table~\ref{table-reentry6d-times} compares this reading with the solution of step 2: the junction times agree within $15$~s (most of them within $5$~s), and $\chi(0)$ within $3^\circ$.
The same structure is found from a constant bank angle of $45^\circ$.
From $90^\circ$, the iteration goes to another local solution, with a flux arc twice as long: like the other methods, this one is local.

\paragraph{Step 2: the multiple-arc TFC method.}
The multiple-arc TFC method has $15$ arcs.
They are the bang arc, the three boundary arcs, the free arcs adjacent to them, which carry the pullback of the neighboring constraint (a bridge between two boundary arcs is split in two halves, each of them carrying one pullback), and the pieces obtained by cutting the long arcs where the smoothed coefficients vary quickly.
These cuts are located around the breakpoints Mach $20$ and Mach $10$ of the tables on the flux and acceleration arcs, and at the middle of the pressure arc.
The smallest degrees are initialized from the solution of step 1 (junction times, $\chi(0)$ and trajectory), and each larger set of degrees from the previous one.

In the structured shooting with events, the unknowns are the switching time of the bang arc, the exit times of the three boundary arcs and $\chi(0)$.
Started from the junction times of the last set, it converges in $30$ evaluations, and its solution (Table~\ref{table-reentry6d-times}) is the reference of the example, with $J_{\mathrm{ref}} = 3.6220471 \times 10^8$~J\,m$^{-2}$.
From the data of step 1, the shooting method does not converge, because the shooting function is not defined there (Section~\ref{discussion}).
A Hermite-Simpson collocation with $N = 800$ intervals, started from the same simulated trajectory, converges to the same structure up to the lateral signs (with sign changes of $\sin \mu$ on the flux and pressure arcs), and to a cost lower by $10^{-4}$ in relative value, the lateral motion being, again, a nearly degenerate direction.
The comparison between the multiple-arc TFC method and shooting is done for a fixed arc structure and fixed branches of the bank sign, and selecting the best lateral strategy is beyond the scope of this example.

\begin{table}[h]
\centering\small
\caption{Reentry problem: junction times (s), final time (s) and initial azimuth (deg) of the structure found by step 1, read from the single-arc solution with $N = 80$, computed by step 2 with $1323$ unknowns, and by the structured shooting started from the step-2 solution.
The angles $\mu_\varphi$, $\mu_a$ and $\mu_P$ are the boundary bank angles of the flux, acceleration and pressure constraints, with the sign of $\sin \mu$.}
\label{table-reentry6d-times}
\setlength{\tabcolsep}{5pt}
\begin{tabular}{llllllllll}\toprule
arc & $\mu = \pi$ & $\mu = 0$ & $\mu_\varphi^-$ & $\mu = 0$ & $\mu_a^+$ & $\mu = 0$ & $\mu_P^+$ & $\mu = 0$ & $\chi(0)$\\\midrule
step 1 & 151.7 & 253.6 & 585.4 & 600.3 & 705.4 & 722.3 & 834.0 & 877.0 & 71.6\\
\midrule
step 2 & 152.92 & 268.43 & 590.38 & 597.84 & 715.14 & 718.57 & 838.64 & 878.74 & 74.70\\
\midrule
shooting & 152.92 & 268.43 & 590.38 & 597.84 & 715.14 & 718.57 & 838.64 & 878.74 & 74.70\\
\bottomrule\end{tabular}
\end{table}

\begin{figure}[h]
\centering
\includegraphics[width=\textwidth]{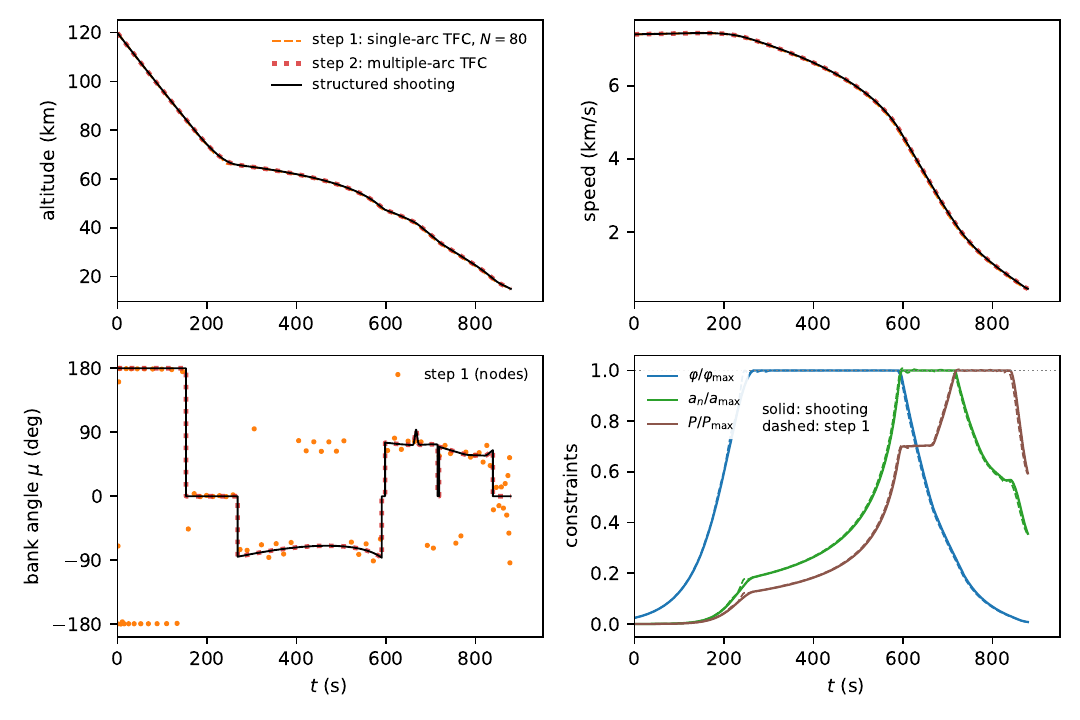}
\caption{Reentry problem: altitude, speed, bank angle and normalized constraints for the single-arc TFC solution of step 1 ($N = 80$, bank angle at the nodes), the multiple-arc TFC solution of step 2 with the largest degrees of Table~\ref{table-reentry6d}, and the structured shooting.}
\label{fig-reentry6d-traj}
\end{figure}

Table~\ref{table-reentry6d} reports the results for the four sets of degrees.
Besides the residual and the junction times, compared with the shooting solution, it reports an off-grid validation: the policy read from the solution ($\chi(0)$, junction times, bank angle on each arc) is integrated without any event, and we measure the terminal errors, the largest constraint violation and the cost error.
The smallest set, started from the rough data of step 1, takes $143$ iterations and three minutes, and the next ones take $14$, $9$ and $9$ iterations, and about fifteen seconds each.
With the largest degrees, the junction times agree with those of the shooting to $6 \times 10^{-8}$~s, for the same smoothed model, the re-integrated trajectory reaches the final altitude to a millimeter, and the cost agrees to $6 \times 10^{-11}$ in relative value.
Of course, this is an agreement between numerical methods on the chosen benchmark, and not an estimate of the uncertainty with respect to the original aerodynamic tables.
Figure~\ref{fig-reentry6d-traj} shows the three solutions and the constraints.

\begin{table}[h]
\centering\footnotesize\setlength{\tabcolsep}{4pt}
\caption{Reentry problem, multiple-arc TFC method for increasing degrees: number of unknowns, residual norm, largest error on the junction times (cuts excluded) and error on the final time with respect to the shooting solution (s), terminal errors on the altitude (m) and on the speed (m/s) and largest constraint violation (normalized) of the re-integrated policy, relative cost error, and solution time (s).}
\label{table-reentry6d}
\begin{tabular}{rrrrrrrrr}\toprule
unknowns & residual & junctions & $t_f$ & $h(t_f) - h_f$ & $v(t_f) - v_f$ & violation & cost & time\\\midrule
647 & $6.0 \times 10^{-5}$ & $7.4 \times 10^{-5}$ & $-6.2 \times 10^{-5}$ & $2.1 \times 10^{-2}$ & $7.6 \times 10^{-4}$ & $2.7 \times 10^{-7}$ & $-4.4 \times 10^{-8}$ & 201\\
877 & $8.0 \times 10^{-7}$ & $5.1 \times 10^{-7}$ & $-1.6 \times 10^{-6}$ & $-4.4 \times 10^{-3}$ & $-1.4 \times 10^{-4}$ & $1.2 \times 10^{-7}$ & $-2.9 \times 10^{-10}$ & 15\\
1109 & $4.1 \times 10^{-8}$ & $2.9 \times 10^{-7}$ & $-8.9 \times 10^{-7}$ & $-2.1 \times 10^{-3}$ & $-7.2 \times 10^{-5}$ & $5.7 \times 10^{-8}$ & $-1.8 \times 10^{-10}$ & 13\\
1323 & $5.4 \times 10^{-9}$ & $5.5 \times 10^{-8}$ & $-7.7 \times 10^{-7}$ & $6.5 \times 10^{-4}$ & $2.2 \times 10^{-5}$ & $1.2 \times 10^{-10}$ & $-6.0 \times 10^{-11}$ & 16\\
\bottomrule\end{tabular}
\end{table}

\subsection{Final remarks}\label{discussion}

Recall that, in the linear-quadratic class, where the dynamics are eliminated, the two TFC methods divide the work as direct and indirect methods do, with the difference that both stages live in the same reduced formulation and that all iterates of both methods are feasible.
The single-arc TFC method locates the boundary arcs approximately, where its slack is small, with an algebraic convergence governed by the junction exponent.
Then the multiple-arc TFC method converts the structure into junction times, with the geometric accuracy of Corollary~\ref{spectral}.
The forced Bryson-Denham variant and the rendezvous problem (whose arcs are not polynomial) illustrate this two-step procedure: a few dozen unknowns provide the accuracy that uniform-mesh collocation reaches with thousands of unknowns, and an upper bound of the optimal cost.
On the unforced Bryson-Denham problem, a Hermite-Simpson transcription whose grid contains the junctions represents the piecewise-cubic solution, and no method can beat it.

For the three examples, Table~\ref{table-shooting} gathers what the shooting method and the multiple-arc TFC method require and what they provide, namely, the number of unknowns, the quantities derived by hand before any code is written, the initial guess, the behavior from initial data far from the solution and the accuracy on the junction times.
On the linear-quadratic problems, the shooting method is the indirect one, on the junction conditions of the maximum principle.
On the reentry problem, the structure fixes the bank angle on every arc, no costate appears, and the shooting method is the structured shooting with events of Section~\ref{reentry}.

\begin{table}[h]
\centering\footnotesize\setlength{\tabcolsep}{4pt}
\caption{Shooting (indirect on the two linear-quadratic problems, structured with events on the reentry problem) against the multiple-arc TFC method: number of unknowns, quantities derived by hand for the problem, initial guess, behavior from initial data far from the solution, accuracy on the junction times; for the multiple-arc TFC method the trials from random junction times are those of the text, in the box of the shooting test.}
\label{table-shooting}
\begin{tabular}{>{\raggedright\arraybackslash}p{2.5cm}>{\raggedright\arraybackslash}p{3.1cm}>{\raggedright\arraybackslash}p{3.1cm}>{\raggedright\arraybackslash}p{4.1cm}}\toprule
 & Bryson-Denham & rendezvous & reentry\\\midrule
\multicolumn{4}{l}{\emph{shooting}}\\
unknowns & $6$ & $5$ & $5$\\
derived by hand & adjoint equations, junction and jump conditions & adjoint equations, jump at the entry & boundary controls, tangency conditions, events\\
initial guess & $p(0)$, $t_1$, $t_2$, $\nu_1$, $\nu_2$ & $p(0)$, $t_1$ & $t_1$, exit times, $\chi(0)$\\
from data far from the solution & $377$ of $500$ random guesses converge & $132$ of $200$ & undefined at the step-1 data; from $1\%$ perturbations converges in $30$ to $62$ evaluations, once out of three to the reference root\\
from the TFC solution & closed form & two Newton steps & $30$ evaluations\\
junction time error & $0$ (closed form) & $0$ (reference, checked in extended precision) & $0$ (reference, integrator tolerance)\\\midrule
\multicolumn{4}{l}{\emph{multiple-arc TFC}}\\
unknowns & $12$ ($N = 6$) & $68$ ($N = 24$) & $1323$\\
derived by hand & none & none & pullbacks, boundary controls\\
initial guess & step 1: threshold on the slack & step 1: threshold on the slack & step 1: hinge fit of the slack\\
from data far from the solution & $44$ of $100$ random junction times converge: all those started on the steep side of the reduced cost, few on its flat side, where the iteration stalls within $0.14$ of the junction and $4 \times 10^{-4}$ of the cost & $35$ of $40$ random entry times and junction values converge & converges from the step-1 data\\
from step 1 & converged in every run & converged in every run & converged\\
admissible iterates & yes & yes & no\\
junction time error & $4 \times 10^{-9}$ & $5 \times 10^{-6}$ & $6 \times 10^{-8}$~s\\
\bottomrule\end{tabular}
\end{table}

The two methods take the arc structure as given, and they differ in what has to be written for each problem.
The indirect shooting method needs the adjoint equations and the junction and jump conditions, derived by hand for each problem.
The resulting code is the smaller one, with six unknowns on the Bryson-Denham problem, whose equations are in closed form, and five on the rendezvous problem.
The multiple-arc TFC method does not involve any adjoint quantity.
Its problem-specific part is the change of coordinates that makes the constraint affine, and the rest (the constrained expressions, the junction values, the residual and its least-squares solution) is the same for every example, the method having more unknowns (twelve and sixty-eight on the same problems).

On the reentry problem, both methods use the boundary controls.
The TFC method adds the pullbacks of the three constraints, and the shooting method adds the tangency conditions and the events that detect the entries, so that the work of derivation is about the same.
Therefore, in the linear-quadratic class, the TFC method is easier to set up, since the constraint is affine as stated and no analysis is needed before writing the code.
Outside this class, the two methods require comparable derivations.

\medskip 

Let us now compare the behavior of the two methods from initial data that are not close to the solution.
The shooting function is defined only where its integration reaches its events, and it has in general several roots.
On the reentry problem, it is not defined at the data read from step 1, although their junction times are within fifteen seconds of the solution.
Indeed, integrated from these data, the policy enters the three boundary arcs with the constraints at $1.01$, $1.15$ and $1.33$ times their bounds (the tangency conditions are the equations to be solved, and they are not yet satisfied), and it goes through the final altitude inside the pressure arc, before its exit time.
Hence the last arc never starts, the function returns a penalty value, and the iteration diverges.
From perturbations of $1\%$ of the solution, it converges in $30$ to $62$ evaluations, but it reaches the reference root in one draw out of three, and another root of the junction equations in the two others, with a final time of $914$~s instead of $879$~s.
On the two linear-quadratic problems, where the arcs are explicit, the shooting function is defined everywhere, and the basins measured in Sections~\ref{bryson} and~\ref{rendezvous} are wide (three random initial guesses out of four, and two out of three), the initial costate being the quantity to be guessed.

For the multiple-arc TFC method, the least-squares residual can be evaluated at the step-1 reentry data, where the event-based shooting fails, but its global convergence is not guaranteed.
On the linear-quadratic problems, the stalls that we observe can be related to the reduced cost, i.e., the cost of the free arcs optimized at fixed junction times.
On the Bryson-Denham problem with $N = 6$, this reduced cost, computed on a $19 \times 19$ grid of the box of the shooting test, has a unique minimum at $(0.3, 0.7)$ and is asymmetric.
When the trial boundary arc is longer than the optimal one, the cost rises quickly (by $4 \times 10^{-2}$ at $(0.25, 0.75)$ and by $6 \times 10^2$ at $(0.05, 0.95)$), because the free arcs are forced to reach the bound too early.
When it is shorter, the cost stays within $10^{-3}$ of the minimum over the whole region $t_1 \geq 0.3$, $t_2 \leq 0.7$, because a free arc then contains a piece of the true boundary arc, which is well approximated by a polynomial of degree $6$.
In this region, the reduced cost has spurious local minima, for instance at $(0.3416, 0.6584)$, $7 \times 10^{-5}$ above the optimum, a point reached from many starts.

From random junction times in the box of the shooting test, each junction time started on its steep side converges ($52$ of $52$ for $t_1$ and $63$ of $63$ for $t_2$), whereas a junction time started on its flat side stalls at a spurious local minimum, unless it starts within $0.02$ of the solution ($8$ of $48$ and $3$ of $37$).
In all, $44$ of $100$ starts converge, and the stalled runs end at junction times off by $0.04$ to $0.14$, with a cost within $4 \times 10^{-4}$ of the optimum.
On the rendezvous problem, with $N = 8$ and with a random entry time and random junction values in the box of the shooting test, $35$ of $40$ starts converge, two stall in the flat valley, within $3 \times 10^{-4}$ of the entry time and $10^{-7}$ of the cost, and three reach other local minima.

The threshold of step 1 overestimates the boundary arc in both linear-quadratic examples, because the slack of the single-arc solution is small on the free arcs near the junctions, where it vanishes with the order given by Proposition~\ref{junction-exponent} (we obtain $[0.2715, 0.7285]$ instead of $[0.3, 0.7]$ for the Bryson-Denham problem, and an entry six percent too early for the rendezvous problem).
This places the initialization of step 2 on the steep side, from which the iteration converged in every run.
On the reentry problem, step 2 converged from the data of step 1.

In these experiments, the robustness of each method depends on the available initialization.
From an initial guess within a few percent of the solution, the shooting method is sharper and cheaper.
Indeed, its convergence is that of Newton's method, and its junction times are exact up to the rounding of the closed form or the tolerance of the integrator, whereas the junction times of the TFC method have errors of $4 \times 10^{-9}$, $5 \times 10^{-6}$ and $6 \times 10^{-8}$~s, their accuracy being limited by the flat valley.
In return, the shooting function has a restricted domain, and one has to choose among its roots.

In the reported tests, the two-step TFC procedure converges without any prior information on the active arcs or on the costate.
Step 1 starts from zero coefficients on the linear-quadratic problems (or from random coefficients: on the Bryson-Denham problem, $100$ of $100$ Gaussian starts with $N = 24$ reach the same cost), and from a trajectory simulated with a constant bank angle on the reentry problem.
These observations concern only the tested initializations.

The two methods are thus complementary: the TFC method provides the structure and an initialization, and the shooting method then refines the junction times.
Started from the costate of the TFC solution on the free arcs (Remark~\ref{dual}) or from its junction times, the shooting method converges in two Newton steps on the rendezvous problem and in thirty evaluations on the reentry problem.
Besides its larger number of unknowns, the TFC method is exact only when the optimal arcs belong to its family (otherwise, it converges geometrically), and each solve is a minimization (or a least-squares problem when the dynamics are not eliminated), and not a root finding.

We have met several limits.
The Gauss-Newton iteration of the single-arc TFC method degenerates once a boundary arc is resolved, since the derivative of the control with respect to the slack vanishes with the slack, and its reduced problem is nonconvex even when the original problem is convex.
The multiple-arc TFC method requires the structure, the junction exponents and an initialization.
It is a local method, whose accuracy on the junction times is limited by the flat valley of the reduced cost, as on the rendezvous problem, where the shooting method started from its solution removes the residual error.
Bounded constraint sets cannot be embedded by independent squared slacks, and the bounds of Section~\ref{ritz} hold only when the dynamics are eliminated.

On the reentry problem, each constraint is exact on the arcs carrying its pullback or its boundary representation.
Started from a trajectory simulated with a constant bank angle, the single-arc TFC method in residual form provides the structure and the junction times to within fifteen seconds, and the multiple-arc TFC method computes the junction times to $10^{-7}$~s with about a thousand unknowns.
But the trial trajectories of both methods are no longer admissible, and the reduced problem is a boundary value problem, and not a minimization.
Moreover, the accuracy on an arc is governed by the regularity of the boundary control, which keeps sharp variations near the breakpoints of the original aerodynamic tables, the iteration of step 1 is local, and the lateral bank sign on the boundary arcs is a nearly degenerate direction, along which the costs of the structures differ by $10^{-4}$ in relative value.

\section{Conclusion}\label{conclusion}

The squared slack turns an affine state constraint into a linear equality, to which TFC applies.
The regularity of the slack function at the junctions, i.e., the H\"older exponent $m/2$ of Proposition~\ref{junction-exponent}, explains the different behaviors of the two methods.
A single polynomial slack converges algebraically, whereas a multiple-arc TFC method, whose free arcs carry the vanishing order of the junction and whose junction times are unknowns, converges geometrically when the arcs are analytic.
In the linear-quadratic class, the identity of Proposition~\ref{ritz-theorem} shows that the cost of every trial trajectory is an upper bound of the optimal cost, and it explains why the square (and not only the tangency) must be kept on the free arcs.

Outside this class, the construction keeps the exactness of the constraints on the arcs adjacent to the boundary arcs, as well as the unknown junction times, and it becomes a boundary value solver for a given structure.
The structure itself (the number and the order of the boundary arcs, and their junction exponents) has to be known, estimated by the single-arc TFC method or by a direct method, or derived from a geometric analysis of the extremals, as in \cite{BonnardFaubourgTrelat2005,Trelat2012JOTA} (see Remark~\ref{structure}).
Cuts inside the arcs are harmless and useful, but a boundary arc that may appear or disappear is outside the present formulation.

We leave for future work the treatment of several constraints that are active simultaneously, the analysis of the residual formulation for general dynamics, and a comparison with adaptive collocation going beyond the rendezvous test reported here.

\paragraph{Acknowledgement.}
E. Tr\'elat was supported by the grant FA8655-25-1-7469 of the European Office of Aerospace Research and Development of the Air Force Office of Scientific.

\small
\bibliographystyle{abbrv}
\bibliography{references}

\begin{thebibliography}{10}

\bibitem{AcikmesePloen2007}
B.~A\c{c}{\i}kme\c{s}e and S.~R. Ploen.
\newblock Convex programming approach to powered descent guidance for {Mars}
  landing.
\newblock {\em Journal of Guidance, Control, and Dynamics}, 30(5):1353--1366,
  2007.

\bibitem{Andersson2019CasADi}
J.~A.~E. Andersson, J.~Gillis, G.~Horn, J.~B. Rawlings, and M.~Diehl.
\newblock {CasADi:} a software framework for nonlinear optimization and optimal
  control.
\newblock {\em Mathematical Programming Computation}, 11(1):1--36, 2019.

\bibitem{Bertsekas1999}
D.~P. Bertsekas.
\newblock {\em Nonlinear Programming}.
\newblock Athena Scientific, Belmont, MA, second edition, 1999.

\bibitem{Betts1998}
J.~T. Betts.
\newblock Survey of numerical methods for trajectory optimization.
\newblock {\em Journal of Guidance, Control, and Dynamics}, 21(2):193--207,
  1998.

\bibitem{Betts2010}
J.~T. Betts.
\newblock {\em Practical Methods for Optimal Control and Estimation Using
  Nonlinear Programming}.
\newblock SIAM, Philadelphia, second edition, 2010.

\bibitem{BonnansHermant2007SICON}
J.~F. Bonnans and A.~Hermant.
\newblock Well-posedness of the shooting algorithm for state constrained
  optimal control problems with a single constraint and control.
\newblock {\em SIAM Journal on Control and Optimization}, 46(4):1398--1430,
  2007.

\bibitem{BonnansHermant2009AIHP}
J.~F. Bonnans and A.~Hermant.
\newblock Second-order analysis for optimal control problems with pure state
  constraints and mixed control-state constraints.
\newblock {\em Annales de l'Institut Henri Poincar\'e C, Analyse Non
  Lin\'eaire}, 26(2):561--598, 2009.

\bibitem{BonnansShapiro2000}
J.~F. Bonnans and A.~Shapiro.
\newblock {\em Perturbation Analysis of Optimization Problems}.
\newblock Springer, New York, 2000.

\bibitem{Bonnard2003Reentry}
B.~Bonnard, L.~Faubourg, G.~Launay, and E.~Tr\'elat.
\newblock Optimal control with state constraints and the space shuttle re-entry
  problem.
\newblock {\em Journal of Dynamical and Control Systems}, 9(2):155--199, 2003.

\bibitem{BonnardFaubourgTrelat2005}
B.~Bonnard, L.~Faubourg, and E.~Tr\'elat.
\newblock Optimal control of the atmospheric arc of a space shuttle and
  numerical simulations with multiple-shooting method.
\newblock {\em Mathematical Models and Methods in Applied Sciences},
  15(1):109--140, 2005.

\bibitem{BonnardTrelat2002COCV}
B.~Bonnard and E.~Tr\'elat.
\newblock Une approche g\'eom\'etrique du contr\^ole optimal de l'arc
  atmosph\'erique de la navette spatiale.
\newblock {\em ESAIM: Control, Optimisation and Calculus of Variations},
  7:179--222, 2002.

\bibitem{BrysonDenhamDreyfus1963}
A.~E. Bryson, W.~F. Denham, and S.~E. Dreyfus.
\newblock Optimal programming problems with inequality constraints. {I}:
  Necessary conditions for extremal solutions.
\newblock {\em AIAA Journal}, 1(11):2544--2550, 1963.

\bibitem{BrysonHo1975}
A.~E. Bryson and Y.-C. Ho.
\newblock {\em Applied Optimal Control: Optimization, Estimation, and Control}.
\newblock Hemisphere Publishing, Washington, DC, revised edition, 1975.

\bibitem{BulirschMontronePesch1991}
R.~Bulirsch, F.~Montrone, and H.~J. Pesch.
\newblock Abort landing in the presence of windshear as a minimax optimal
  control problem, part 2: Multiple shooting and homotopy.
\newblock {\em Journal of Optimization Theory and Applications},
  70(2):223--254, 1991.

\bibitem{ByczkowskiRao2024}
C.~A. Byczkowski and A.~V. Rao.
\newblock Method for solving state-path constrained optimal control problems
  using adaptive {R}adau collocation.
\newblock {\em Optimal Control Applications and Methods}, 45(3):1199--1222,
  2024.

\bibitem{CaillauCotsGergaud2012}
J.-B. Caillau, O.~Cots, and J.~Gergaud.
\newblock Differential continuation for regular optimal control problems.
\newblock {\em Optimization Methods and Software}, 27(2):177--196, 2012.

\bibitem{CaillauCotsMartinon2022ct}
J.-B. Caillau, O.~Cots, and P.~Martinon.
\newblock ct: control toolbox -- numerical tools and examples in optimal
  control.
\newblock {\em IFAC-PapersOnLine}, 55(16):13--18, 2022.

\bibitem{Clarke2013}
F.~Clarke.
\newblock {\em Functional Analysis, Calculus of Variations and Optimal
  Control}, volume 264 of {\em Graduate Texts in Mathematics}.
\newblock Springer, London, 2013.

\bibitem{ClohessyWiltshire1960}
W.~H. Clohessy and R.~S. Wiltshire.
\newblock Terminal guidance system for satellite rendezvous.
\newblock {\em Journal of the Aerospace Sciences}, 27(9):653--658, 1960.

\bibitem{DAmbrosio2022Landing}
A.~D'Ambrosio, E.~Schiassi, H.~Johnston, F.~Curti, D.~Mortari, and R.~Furfaro.
\newblock Time-energy optimal landing on planetary bodies via theory of
  functional connections.
\newblock {\em Advances in Space Research}, 69(12):4198--4220, 2022.

\bibitem{deAlmeida2026Segmentation}
A.~de~Almeida.
\newblock Segmentation of the spacecraft transfer problem based on
  overdetermined and continuity constraints using the theory of functional
  connections.
\newblock {\em Acta Astronautica}, 242:251--263, 2026.

\bibitem{DeVoreLorentz1993}
R.~A. DeVore and G.~G. Lorentz.
\newblock {\em Constructive Approximation}, volume 303 of {\em Grundlehren der
  mathematischen Wissenschaften}.
\newblock Springer, Berlin, 1993.

\bibitem{DingWright2025}
L.~Ding and S.~J. Wright.
\newblock On squared-variable formulations.
\newblock {\em SIAM Journal on Optimization}, 35(4):2265--2293, 2025.

\bibitem{DrozdFurfaro2025}
K.~Drozd and R.~Furfaro.
\newblock Extreme theory of functional connections with receding horizon
  control for aerospace applications.
\newblock {\em Mathematics}, 13(22):3717, 2025.

\bibitem{DrozdFurfaroDAmbrosio2024}
K.~Drozd, R.~Furfaro, and A.~D'Ambrosio.
\newblock A theory of functional connections-based {\emph{hp}}-adaptive mesh
  refinement algorithm for solving hypersensitive two-point boundary-value
  problems.
\newblock {\em Mathematics}, 12(9):1360, 2024.

\bibitem{Drozd2021}
K.~Drozd, R.~Furfaro, E.~Schiassi, H.~Johnston, and D.~Mortari.
\newblock Energy-optimal trajectory problems in relative motion solved via
  theory of functional connections.
\newblock {\em Acta Astronautica}, 182:361--382, 2021.

\bibitem{EichmeirNachbagauerSteiner2025}
P.~Eichmeir, K.~Nachbagauer, and W.~Steiner.
\newblock The use of slack variables in the adjoint method handling inequality
  constraints in optimal control and the application to tumor drug dosage.
\newblock {\em Journal of Computational and Nonlinear Dynamics}, 20(2):021004,
  2025.

\bibitem{Fliess1995}
M.~Fliess, J.~L\'evine, P.~Martin, and P.~Rouchon.
\newblock Flatness and defect of non-linear systems: introductory theory and
  examples.
\newblock {\em International Journal of Control}, 61(6):1327--1361, 1995.

\bibitem{FourerGayKernighan2003}
R.~Fourer, D.~M. Gay, and B.~W. Kernighan.
\newblock {\em {AMPL}: A Modeling Language for Mathematical Programming}.
\newblock Duxbury Press / Brooks/Cole, Pacific Grove, CA, second edition, 2003.

\bibitem{FukudaFukushima2017}
E.~H. Fukuda and M.~Fukushima.
\newblock A note on the squared slack variables technique for nonlinear
  optimization.
\newblock {\em Journal of the Operations Research Society of Japan},
  60(3):262--270, 2017.

\bibitem{FurfaroMortari2020}
R.~Furfaro and D.~Mortari.
\newblock Least-squares solution of a class of optimal space guidance problems
  via theory of connections.
\newblock {\em Acta Astronautica}, 168:92--103, 2020.

\bibitem{GraichenPetit2009}
K.~Graichen and N.~Petit.
\newblock Incorporating a class of constraints into the dynamics of optimal
  control problems.
\newblock {\em Optimal Control Applications and Methods}, 30(6):537--561, 2009.

\bibitem{HaberkornMartinonGergaud2004}
T.~Haberkorn, P.~Martinon, and J.~Gergaud.
\newblock Low thrust minimum-fuel orbital transfer: A homotopic approach.
\newblock {\em Journal of Guidance, Control, and Dynamics}, 27(6):1046--1060,
  2004.

\bibitem{HargravesParis1987}
C.~Hargraves and S.~Paris.
\newblock Direct trajectory optimization using nonlinear programming and
  collocation.
\newblock {\em Journal of Guidance, Control, and Dynamics}, 10(4):338--342,
  1987.

\bibitem{HartlSethiVickson1995}
R.~F. Hartl, S.~P. Sethi, and R.~G. Vickson.
\newblock A survey of the maximum principles for optimal control problems with
  state constraints.
\newblock {\em SIAM Review}, 37(2):181--218, 1995.

\bibitem{JacobsonLele1969}
D.~H. Jacobson and M.~M. Lele.
\newblock A transformation technique for optimal control problems with a state
  variable inequality constraint.
\newblock {\em IEEE Transactions on Automatic Control}, 14(5):457--464, 1969.

\bibitem{JacobsonLeleSpeyer1971}
D.~H. Jacobson, M.~M. Lele, and J.~L. Speyer.
\newblock New necessary conditions of optimality for control problems with
  state-variable inequality constraints.
\newblock {\em Journal of Mathematical Analysis and Applications},
  35(2):255--284, 1971.

\bibitem{JohnstonLeakeMortari2020Ineq}
H.~Johnston, C.~Leake, and D.~Mortari.
\newblock An analysis of the theory of functional connections subject to
  inequality constraints.
\newblock In {\em Advances in the Astronautical Sciences}, volume 171, pages
  325--337, 2020.
\newblock Paper AAS 19-732.

\bibitem{JohnstonMortari2021Hybrid}
H.~Johnston and D.~Mortari.
\newblock Least-squares solutions of boundary-value problems in hybrid systems.
\newblock {\em Journal of Computational and Applied Mathematics}, 393:113524,
  2021.

\bibitem{Johnston2020Descent}
H.~Johnston, E.~Schiassi, R.~Furfaro, and D.~Mortari.
\newblock {Fuel-Efficient} powered descent guidance on large planetary bodies
  via theory of functional connections.
\newblock {\em The Journal of the Astronautical Sciences}, 67(4):1521--1552,
  2020.

\bibitem{Kelly2017}
M.~Kelly.
\newblock An introduction to trajectory optimization: How to do your own direct
  collocation.
\newblock {\em SIAM Review}, 59(4):849--904, 2017.

\bibitem{LeakeJohnstonMortari2022}
C.~Leake, H.~Johnston, and D.~Mortari.
\newblock {\em The Theory of Functional Connections: A Functional Interpolation
  Framework with Applications}.
\newblock Lulu, Morrisville, NC, 2022.

\bibitem{LeakeMortari2020Deep}
C.~Leake and D.~Mortari.
\newblock Deep theory of functional connections: A new method for estimating
  the solutions of partial differential equations.
\newblock {\em Machine Learning and Knowledge Extraction}, 2(1):37--55, 2020.

\bibitem{LeakeThesis2021}
C.~D. Leake.
\newblock {\em The Multivariate Theory of Functional Connections: An
  $n$-Dimensional Constraint Embedding Technique Applied to Partial
  Differential Equations}.
\newblock PhD thesis, Texas A\&M University, College Station, TX, 2021.

\bibitem{LuLiu2013}
P.~Lu and X.~Liu.
\newblock Autonomous trajectory planning for rendezvous and proximity
  operations by conic optimization.
\newblock {\em Journal of Guidance, Control, and Dynamics}, 36(2):375--389,
  2013.

\bibitem{MaiMortari2022}
T.~Mai and D.~Mortari.
\newblock Theory of functional connections applied to quadratic and nonlinear
  programming under equality constraints.
\newblock {\em Journal of Computational and Applied Mathematics}, 406:113912,
  2022.

\bibitem{Malisani2023SICON}
P.~Malisani.
\newblock Interior point methods in optimal control problems of affine systems:
  Convergence results and solving algorithms.
\newblock {\em SIAM Journal on Control and Optimization}, 61(6):3390--3414,
  2023.

\bibitem{Malisani2024COCV}
P.~Malisani.
\newblock Interior point methods in optimal control.
\newblock {\em ESAIM: Control, Optimisation and Calculus of Variations}, 30:59,
  2024.

\bibitem{Maurer1977}
H.~Maurer.
\newblock On optimal control problems with bounded state variables and control
  appearing linearly.
\newblock {\em SIAM Journal on Control and Optimization}, 15(3):345--362, 1977.

\bibitem{MaurerGillessen1975}
H.~Maurer and W.~Gillessen.
\newblock Application of multiple shooting to the numerical solution of optimal
  control problems with bounded state variables.
\newblock {\em Computing}, 15(2):105--126, 1975.

\bibitem{MilamMushambiMurray2000}
M.~B. Milam, K.~Mushambi, and R.~M. Murray.
\newblock A new computational approach to real-time trajectory generation for
  constrained mechanical systems.
\newblock In {\em Proceedings of the 39th IEEE Conference on Decision and
  Control}, volume~1, pages 845--851, 2000.

\bibitem{Mortari2017LS}
D.~Mortari.
\newblock {Least-Squares} solution of linear differential equations.
\newblock {\em Mathematics}, 5(4):48, 2017.

\bibitem{Mortari2017Connections}
D.~Mortari.
\newblock The theory of connections: Connecting points.
\newblock {\em Mathematics}, 5(4):57, 2017.

\bibitem{Mortari2025MinParam}
D.~Mortari.
\newblock Leveraging minimum parameterization in the theory of functional
  connections.
\newblock {\em Lobachevskii Journal of Mathematics}, 46(8):4106--4123, 2025.

\bibitem{MortariJohnstonSmith2019}
D.~Mortari, H.~Johnston, and L.~Smith.
\newblock High accuracy least-squares solutions of nonlinear differential
  equations.
\newblock {\em Journal of Computational and Applied Mathematics}, 352:293--307,
  2019.

\bibitem{MortariLeake2019}
D.~Mortari and C.~Leake.
\newblock The multivariate theory of connections.
\newblock {\em Mathematics}, 7(3):296, 2019.

\bibitem{Pesch1994}
H.~J. Pesch.
\newblock A practical guide to the solution of real-life optimal control
  problems.
\newblock {\em Control and Cybernetics}, 23(1--2):7--60, 1994.

\bibitem{Schiassi2021XTFC}
E.~Schiassi, R.~Furfaro, C.~Leake, M.~De~Florio, H.~Johnston, and D.~Mortari.
\newblock Extreme theory of functional connections: A fast physics-informed
  neural network method for solving ordinary and partial differential
  equations.
\newblock {\em Neurocomputing}, 457:334--356, 2021.

\bibitem{Taylor1973}
J.~G. Taylor.
\newblock A squared-variable transformation approach to nonlinear programming
  optimality conditions.
\newblock {\em Naval Research Logistics Quarterly}, 20(1):25--39, 1973.

\bibitem{Trefethen2019ATAP}
L.~N. Trefethen.
\newblock {\em Approximation Theory and Approximation Practice, Extended
  Edition}.
\newblock SIAM, Philadelphia, 2019.

\bibitem{Trelat2003DCDS}
E.~Tr\'elat.
\newblock Optimal control of a space shuttle, and numerical simulations.
\newblock {\em Discrete and Continuous Dynamical Systems}, suppl.:839--848,
  2003.
\newblock Dynamical systems and differential equations (Wilmington, NC, 2002).

\bibitem{Trelat2012JOTA}
E.~Tr\'elat.
\newblock Optimal control and applications to aerospace: Some results and
  challenges.
\newblock {\em Journal of Optimization Theory and Applications},
  154(3):713--758, 2012.

\bibitem{Valentine1937}
F.~A. Valentine.
\newblock The problem of {Lagrange} with differential inequalities as added
  side conditions.
\newblock In {\em Contributions to the Calculus of Variations 1933--1937},
  pages 407--448. University of Chicago Press, Chicago, 1937.
\newblock Reprinted in: Traces and Emergence of Nonlinear Programming,
  Birkh\"auser, Basel, 2014, pp. 313--357, doi:10.1007/978-3-0348-0439-4\_16.

\bibitem{VanLoock2015}
W.~Van~Loock, G.~Pipeleers, and J.~Swevers.
\newblock B-spline parameterized optimal motion trajectories for robotic
  systems with guaranteed constraint satisfaction.
\newblock {\em Mechanical Sciences}, 6(2):163--171, 2015.

\bibitem{Vinter2000}
R.~Vinter.
\newblock {\em Optimal Control}.
\newblock Systems \& Control: Foundations \& Applications. Birkh\"auser,
  Boston, 2000.

\bibitem{Virtanen2020SciPy}
P.~Virtanen et~al.
\newblock {SciPy} 1.0: fundamental algorithms for scientific computing in
  python.
\newblock {\em Nature Methods}, 17(3):261--272, 2020.

\bibitem{WachterBiegler2006}
A.~W\"achter and L.~T. Biegler.
\newblock On the implementation of an interior-point filter line-search
  algorithm for large-scale nonlinear programming.
\newblock {\em Mathematical Programming}, 106(1):25--57, 2006.

\end{thebibliography}

\end{document}